\documentclass[english,opre,nonblindrev]{informs3_arxiv}
\usepackage[T1]{fontenc}
\usepackage[latin9]{inputenc}
\usepackage{geometry}
\usepackage{xcolor}
\usepackage{array}
\usepackage{float}
\usepackage{booktabs}
\usepackage{bm}
\usepackage{amsmath}
\usepackage{amssymb}
\usepackage{graphicx}
\usepackage{esint}
\usepackage[authoryear]{natbib}
\usepackage[hidelinks]{hyperref}

\usepackage{enumitem}

\definecolor{darkbrown}{rgb}{0.4, 0.26, 0.13}
\definecolor{darkcandyapplered}{rgb}{0.64, 0.0, 0.0}
\definecolor{darkolivegreen}{rgb}{0.33, 0.42, 0.18}
\definecolor{darkred}{rgb}{0.55, 0.0, 0.0}
\definecolor{davysgrey}{rgb}{0.33, 0.33, 0.33}
\definecolor{darkgreen}{rgb}{0.0, 0.2, 0.13}
\definecolor{darkorchid}{rgb}{0.6, 0.2, 0.8}
\definecolor{darkblue}{rgb}{0.0, 0.0, 0.55}

\newcommand\mycolor{black}

\makeatletter

\providecommand{\tabularnewline}{\\}

\DoubleSpacedXI 

\usepackage{endnotes}
\let\footnote=\endnote

\usepackage{pgfplots}

\usepackage{natbib}
 \bibpunct[, ]{(}{)}{,}{a}{}{,}%
 \def\bibfont{\small}%
\TheoremsNumberedThrough     

\newcommand{\E}{\mathbb{E}}

\newcommand{\R}{\mathbb{R}}
\newcommand{\I}{\mathbbm{1}}

\newcommand{\map}{\rightarrow}
\newcommand{\CC}{\mathbb{C}}

\newcommand{\fx}{f_{\mathbf{X}}}
\newcommand{\fy}{f_{\mathbf{Y}}}

\newcommand{\fv}{f_{\mathbf{V}}}
\newcommand{\FX}{F_{\mathbf{X}}}
\newcommand{\FY}{F_{\mathbf{Y}}}

\newcommand{\FV}{F_{\mathbf{V}}}
\newcommand{\pX}{P_{\mathbf{X}}}
\newcommand{\pY}{P_{\mathbf{Y}}}

\newcommand{\pV}{P_{\mathbf{V}}}

\newcommand{\yy}{{\mathbf{y}}}
\newcommand{\zz}{{\mathbf{z}}}

\newcommand{\YY}{{\mathbf{Y}}}

\newcommand{\1}{\mathbbm{1}}

\usepackage{todonotes}

\usepackage[none]{hyphenat}
\makeatletter

\usepackage{endnotes}
\let\footnote=\endnote

\usepackage{natbib}
 \bibpunct[, ]{(}{)}{,}{a}{}{,}%
 \def\bibfont{\small}%
\usepackage{amsmath,bbm}
\usepackage{footnote}

\def\footnoterule{\kern-3\p@
  \hrule \@width 6.5in \kern 2.6\p@} 

\makeatother

\usepackage{babel}
\begin{document}
\title{Scores for Multivariate Distributions and Level Sets}

\RUNAUTHOR{Meng, Taylor, Ben Taieb and Li}
\RUNTITLE{Scoring Function, Distribution, Level Set}

\ARTICLEAUTHORS{%
\AUTHOR{Xiaochun Meng}
\AFF{University of Sussex Business School, University of Sussex, UK, \EMAIL{xiaochun.meng@sussex.ac.uk}}
\AUTHOR{ James W. Taylor}
\AFF{Sa\"{i}d Business School, University of Oxford, UK, \EMAIL{james.taylor@sussex.ac.uk}}
\AUTHOR{ Souhaib Ben Taieb}
\AFF{Department of Computer Science, University of Mons, Belgium, \EMAIL{souhaib.bentaieb@umons.ac.be}}
\AUTHOR{  Siran Li\footnotemark}
\AFF{School of Mathematical Sciences $\&$ IMA-Shanghai, Shanghai Jiao Tong University, Shanghai, China, \EMAIL{sl4025@nyu.edu}}
\footnotetext{\linespread{0.50}\selectfont{}The research of Siran Li is supported by NSFC (National Natural Science Foundation of China) Grant No.~12201399 and Shanghai Frontier Research Institute for Modern Analysis.}
}

\abstract{Forecasts of multivariate probability distributions are required for a variety of applications. Scoring rules enable the evaluation of forecast accuracy, and comparison between forecasting methods. We propose a theoretical framework for scoring rules for multivariate distributions, which encompasses the existing quadratic score and multivariate continuous ranked probability score. We demonstrate how this framework can be used to generate new scoring rules. In some multivariate contexts, it is a forecast of a level set that is needed, such as a density level set for anomaly detection or the level set of the cumulative distribution as a measure of risk. This motivates consideration of scoring functions for such level sets. For univariate distributions, it is well-established that the continuous ranked probability score can be expressed as the integral over a quantile score. We show that, in a similar way, scoring rules for multivariate distributions can be decomposed to obtain scoring functions for level sets. Using this, we present scoring functions for different types of level set, including density level sets and level sets for cumulative distributions. To compute the scores, we propose a simple numerical algorithm. {\color{\mycolor}We perform a simulation study to support our proposals, and we use real data to illustrate usefulness for forecast combining and CoVaR estimation.}}


\maketitle


\section{Introduction}

Forecasts should be probabilistic, as this allows forecasters to communicate
their uncertainty and hence better support decision-making. In many
applications, a forecast of a multivariate probability distribution
is needed. For example, \citet{berrocal2010spatio} consider the impact
of extreme weather risk on road maintenance by estimating the joint
distribution of temperature and precipitation to enable the prediction
of ice formation on road surfaces. \citet{danaher2011modeling} show
how accurate estimation of the joint distribution of website visit
duration and expenditure can lead to improved distributional forecasts
for the latter. \citet{jeon2012using} forecast the joint distribution
of wind velocity in perpendicular directions, which they use to generate
distributional forecasts of wind power. \citet{diks2014comparing}
model the joint distribution of sets of exchange rates to capture
tail dependence and its impact on common extreme appreciation and
depreciation of the currencies.

To support such applications, \citet{gneiting2014probabilistic} describe \textcolor{\mycolor}{the} need for decision-theoretically principled methods
for evaluating forecasts of multivariate distributions. Distributional
forecast accuracy should be evaluated by maximizing sharpness subject
to calibration. Sharpness relates to the concentration of the probabilistic
forecast, while calibration concerns its statistical consistency with
the data. A score summarizes both calibration and sharpness, and can
be used to compare forecasts from competing methods. In addition,
a score can be used as the objective function in model estimation
(\citealt{jose2017percentage}).

The continuous ranked probability score (CRPS) is widely used for
univariate distributions (\citealt{grushka2017quantile}). For multivariate
distributions, the energy score of \citet{gneiting2007strictly} has
become \textcolor{\mycolor}{popular} as a multivariate generalization
of the CRPS. An alternative generalization is the multivariate CRPS
(MCRPS), introduced by \citet{gneiting2007strictly}. This score has
only been considered further by \citet{yuen2014crps}, who use it
for estimation in extreme value theory. Another proposal is the quadratic
score, which can be used for univariate and multivariate densities,
but is distinct from the log score (see \citealt{gneiting2007strictly}).

Often the object of interest of a probability distribution is a \textcolor{\mycolor}{statistic, or a property}, such as a quantile of a univariate distribution. Scores for
quantiles have been thoroughly studied (see, for example, \citealt{grushka2017quantile,jose2009evaluating}).
However, there is a far less developed literature on scores for level
sets, which can be viewed as multivariate generalizations of quantiles. {\color{\mycolor} Roughly speaking, for a function $g$, an $\alpha$
level set contains all the points in $\mathbb{R}^{d}$ where $g$
is greater than or equal to $\alpha$. For example, if $g$ is a cumulative distribution function (CDF), the $\alpha$ level set of $g$ is called an
$\alpha$ CDF level set.} To the best of our knowledge, the only score
available for any form of multivariate level set is the excess mass
score for density level sets (see, for example, \citealt{hartigan1987estimation}).
Estimates of such level sets are important for a variety of applications. 

In the context of financial and insurance risk management, \citet{embrechts2006bounds}
propose the boundary of a \textcolor{\mycolor}{CDF} level
set as a multivariate extension of value at risk (VaR). VaR helps
institutions to decide the quantity of assets required to cover potential
losses. The alternative risk measure of \citet{cousin2013multivariate}
is the expectation of the underlying portfolio of risks, conditional
on it lying on the boundary of the level set. 

CDF level sets have been used frequently in environmental risk management,
particularly in relation to hydrology and coastal management. For
example, \citet{corbella2012multivariate} use CDF level sets of the
joint distribution of storm duration and wave height to quantify risk
in relation to coastal erosion. For similar data, \citet{salvadori2014practical}
show how the CDF level set can be used to find critical values for
storm duration and wave height, which are then used to decide the
appropriate weight of concrete units in the design of a breakwater.
In a bivariate flood hazard analysis, \citet{moftakhari2017compounding}
use the boundaries of CDF level sets to assess risks associated with
the compounding effects of river flooding and sea level rise. Each
boundary consists of pairs of values of coastal water level and fluvial
flow for which there is a common probability of exceedance by at least
one of the two variables. This is viewed as a failure probability,
which can be used to classify future hazards in a warming climate.
CDF level sets are also \textcolor{\mycolor}{central to} the four definitions
of \citet{salvadori2016multivariate} for different forms of hazard
scenarios in extreme natural phenomena. \citet{newman2017review}
review decision support systems targeted at reducing natural hazard
risk. CDF level sets have also been used in the judgemental estimation
of bivariate distributions (\citealt{abbas2010assessing}). The joint
distribution is constructed from CDF level sets judgmentally estimated
with the aid of a sequence of binary choices.

Density level sets, which are \textcolor{\mycolor}{associated with} probability density functions
(PDFs), have been used for anomaly detection and cluster analysis
(\citealt{chen2017density}). \citet{steinwart2005classification}
explain that a natural way to identify anomalies is to use density
level sets because they precisely describe the concentration of the
probability distribution. An observation is deemed to be an anomaly
if it lies outside the level set corresponding to a chosen value of
the density, which serves as the anomaly tolerance level (\citealt{rigollet2009optimal}).
Another common use for density level sets is cluster analysis, which
is used extensively for preliminary data analysis, classification
and data reduction (\citealt{menardi2016review}). For a chosen value
of the density, the observations that lie within the same connected
component of the level set can be considered to be homogeneous (see
\citealt{hartigan1975clustering,rinaldo2010generalized,rinaldo2012stability}).
This form of clustering is closely related to modal clustering (see
\citealt{menardi2016review}), and indeed density level sets have
also been used as the basis for tests of multimodality (see \citealt{muller1991excess}).

In this paper, we aim to extend the literature on the evaluation of
forecasts of multivariate distributions and their level sets. For
a univariate distribution, the CRPS can be expressed as an integral
over a quantile score (\citealt{laio2007verification}). We generalize
this to the multivariate context. In doing so, we make a number of
contributions. First, we propose a theoretical framework that links
the quadratic score, CRPS and MCRPS. This framework can be used to
generate new scores for multivariate distributions, and we demonstrate
this by developing a score based on lower partial moments (LPMs) (\citealt{price1982variance,briec2010portfolio,anthonisz2012asset}).
Second, we show that by decomposing the quadratic score, MCRPS, and
our new score, we obtain new scores for \textcolor{\mycolor}{density, CDF and LPM level sets}, respectively. These scores encompass the existing scores for
density level sets, and the full class of quantile scores. Finally,
to compute the \textcolor{\mycolor}{different} scores, we propose a simulation-based numerical
approach.

\begin{figure}[H]
\centering{}\caption{\linespread{0.50}\selectfont{} 3-D illustration of the density, CDF
and LPM level sets for the bivariate normal distribution
with zero means, unit variances, and covariance 0.5.\label{fig:Illustration level sets 3D}}
\includegraphics[scale=0.2]{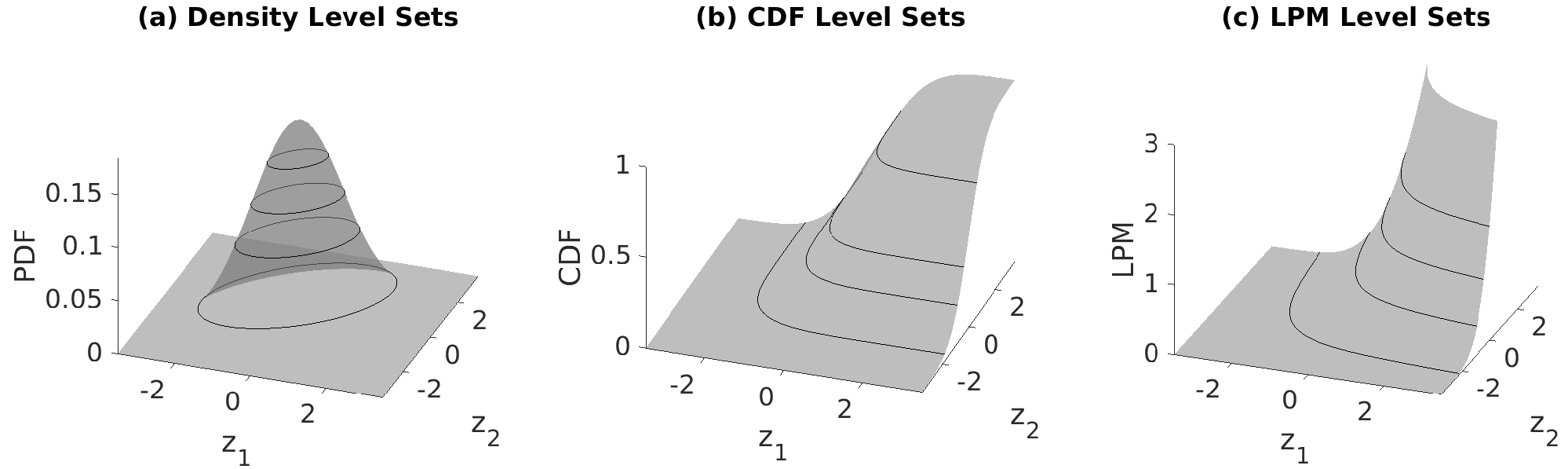}
\end{figure}

To accompany our development of the level set scores, we use a running
illustrative example based on a bivariate normal distribution with
zero means, unit variances, and covariance 0.5. We refer to this as
the data generating process (DGP). For the DGP, Figure~\ref{fig:Illustration level sets 3D}
presents 3-D illustrative examples of density, CDF and LPM level sets.
These level sets will be formally defined in Section~\ref{sec:Scoring Functions for Multivariate Level Sets}.

We note that scores for vector-valued properties have been extensively
studied (see, for example, \citealt{abernethy2012characterization,frongillo2015vector,lambert2019elicitation}).
Examples of vector-valued properties include the mean and covariance.
That literature is distinct from our consideration of level sets,
because level sets are set-valued, rather than vector-valued: {\color{\mycolor}a level
set can be any Borel set, e.g., an ellipsoid or hyperspace, which
cannot be described by a vector.}

{\color{\mycolor} Section~$\ref{sec:Preliminaries}$ describes notations and conventions used in this paper. In Section~$\ref{sec:A-Framework-for scoring functions}$, we present our new theoretical framework for building scoring rules for multivariate distributions, and we demonstrate how this framework can be used to generate new scoring rules. Section~$\ref{sec:Scoring Functions for Multivariate Level Sets}$ shows that the scores presented in Section~$\ref{sec:A-Framework-for scoring functions}$ can be decomposed to obtain scores for different types of level sets. In Section~$\ref{sec:Simulation-Study}$, after discussing how to compute the scores, we provide support for the new scores using simulated data. Section 6 provides two practical applications of the new scores. First, we use our new score for distributions to estimate weights in a combination of distributional forecasts taken from the ECB Survey of Professional Forecasters. The new score is particularly convenient for this application as it enables the combining weights to be estimated using a quadratic optimization algorithm, which brings both theoretical and computational advantages. Second, using financial returns data, we show how CDF level sets computed with our scores provide better estimates of conditional value at risk (CoVaR), a widely-used measure of systemic risk. Section~$\ref{sec:Conclusion}$ concludes the paper.}

\section{Preliminaries\label{sec:Preliminaries}}

In this section, we explain some notations and conventions used in
the subsequent parts of the paper. We identify a vector in $\R^{d}$
with a $d\times1$ (column) matrix. The boldface lower-case letters
$\mathbf{z},\mathbf{s},\mathbf{t},\ldots$ designate points in $\R^{d}$;
we use $z_{j}$ to denote each coordinate; and $\|\bullet\|$ is the
Euclidean modulus. For $\mathbf{z}\in\R$ or $\CC$ we write $\|\mathbf{z}\|\equiv|\mathbf{z}|$.
Given any Borel set $A\subseteq\R^{d}$, $\partial A$ denotes its
topological boundary. Let $\lambda$ be a Borel measure on $\mathbb{R}^{d}$.
We say that a function $g:\mathbb{R}^{d}\rightarrow\mathbb{R}^{m}$
is integrable with respect to $\lambda$ if $\int_{\mathbb{R}^{d}}\|g({\bf z})\|\,d\lambda({\bf z})<\infty$.
Let $\delta_{{\bf y}}$ be the Dirac delta measure with the mass at
$\mathbf{y}$, and let $\mathcal{L}^{d}$ be the $d$-dimensional
Lebesgue measure. The $L^{2}$ norm of a function $g$ with respect
to the measure $\lambda$ is defined as $\mathbf{\int}_{\mathbb{R}^{d}}\|g({\bf z})\|^{2}\,d\lambda({\bf z})$.
For two univariate functions $g_{1}$ and $g_{2}$, i.e. $m=1$, the
convolution of $g_{1}$ and $g_{2}$ is defined as $\left(g_{1}*g_{2}\right)(\mathbf{z})=\mathbf{\int}_{\R^{d}}g_{1}(\mathbf{s})g_{2}(\mathbf{z}-\mathbf{s})d\mathbf{s}$.

We use capital letters $X,Y,V,\ldots$ for univariate random variables,
and boldface capital letters $\mathbf{X},\mathbf{Y},\mathbf{V},\ldots$
for multivariate random variables. As in the literature (see, for
example, \citealt{gneiting2007strictly}), we denote by $P_{{\bf X}},P_{{\bf Y}},P_{{\bf V}}\ldots$
the probability measures of $\mathbf{X},\mathbf{Y},\mathbf{V},\ldots$,
and by $\mathcal{V}^{d}$ a convex class of probability measures.
For distributions $\pX,\pY,\pV,\ldots$, denote by $\FX,\FY,\FV,\ldots$
their CDFs, and denote by $\fx,\fy,\fv,\ldots$ their {\color{\mycolor} generalized probability
density functions}.

For continuous distributions, $\fx,\fy,\fv,\ldots$
are the PDFs; for discrete distributions, $\fx,\fy,\fv,\ldots$ are
the probability mass functions; and $\fx,\fy,\fv,\ldots$ can also
represent the probability measures for linear combinations of continuous
and discrete distributions. For simplicity, here and hereafter, we
refer to $\fx,\fy,\fv,\ldots$ simply as the PDFs. We write $q_{Y}(\alpha)=\inf\{z\in\mathbb{R}|F_{Y}\geq\alpha\}$
for the $\alpha$ quantile of a univariate random variable $Y$ for
$\alpha\in[0,1]$.

\section{A New Framework for Scoring Rules\label{sec:A-Framework-for scoring functions}}

Let $P_{\mathbf{Y}}\in\mathcal{V}^{d}$ be the unknown distribution
of a random variable $\mathbf{Y}$ we want to study. A key problem
is to evaluate the accuracy of $P_{\mathbf{X}}$, an estimate of $P_{\mathbf{Y}},$
given a finite collection of realizations of $\mathbf{Y}$, labeled
as $\{\mathbf{y}^{i}\}_{i=1,2,...T}$. A central tool developed for
this purpose is the \textit{scoring rule}:
\[
S(P_{\mathbf{X}},\mathbf{y}):\mathcal{V}^{d}\times\R^{d}\longrightarrow\R.
\]
It is said to be \textit{proper} if $\E_{P_{\mathbf{Y}}}\left[S(P_{\mathbf{X}},\bullet)\right]$,
the expectation of $S(P_{\mathbf{X}},\bullet)$ with respect to $P_{\mathbf{Y}}$,
is minimized by $P_{\mathbf{Y}}$, and \textit{strictly proper} if
$P_{\mathbf{Y}}$ is the unique minimizer (\citealt{gneiting2007strictly}).

{\color{\mycolor} It can be intuitive to understand a score $S$ through
its divergence, which is defined as $\E_{P_{\mathbf{Y}}}\left[S(P_{\mathbf{X}},\bullet)\right]-\E_{P_{\mathbf{Y}}}\left[S(P_{\mathbf{Y}},\bullet)\right]$.
The divergence is always non-negative, and is equal to 0 if $P_{\mathbf{X}}=P_{\mathbf{Y}}$.
Hence, it essentially defines a measure of ``distance'' between distributions.
Scores that have the same divergence lead to the same measure of distance,
so are often viewed as \textit{equivalent}.}

In the next section, we propose a class of $L^{2}$ scoring rules
for distributions. In Section~\ref{subsec: Examples of Scoring Functions},
we show that the new class of scoring rules encompasses the quadratic
score, CRPS and MCRPS. We also show that this framework can easily
be used to generate other scoring rules, and we demonstrate this by
proposing a scoring rule based on lower partial moments.

\subsection{A Novel Method for Constructing Scoring Rules\label{subsec: scoring fn}}

{\color{\mycolor}
The quadratic score is a well-known proper scoring rule, defined as $\int_{\R^{d}}f_{\mathbf{X}}^{2}(\mathbf{z})\,d{\bf z}-2f_{\mathbf{X}}(\mathbf{y})$, whose divergence function is the $L^{2}$ distance between PDFs, i.e. $\int_{\R^{d}}\left(f_{\mathbf{X}}(\mathbf{z})-f_{\mathbf{Y}}(\mathbf{z})\right)^{2}d{\bf z}$. 
Despite its popularity, the quadratic score has the notable weakness that sometimes it cannot distinguish the relative accuracy
between misspecified estimates. For example, if $f_{\mathbf{X}}$ and $f_{\mathbf{Y}}$ have disjoint support, then the divergence is simply $\int_{\R^{d}}f^2_{\mathbf{X}}(\mathbf{z})d{\bf z}+\int_{\R^{d}}f^2_{\mathbf{Y}}(\mathbf{z})d{\bf z}$. Consider the case where $f_{\mathbf{Y}}{(\mathbf{z})}=\delta_{\mathbf{0}}(\mathbf{z})$ representing all the probability mass at 0, and  $f_{\mathbf{X}}(\mathbf{z})=\delta_{\mathbf{x}_0}(\mathbf{z})$ representing all the probability mass at $\mathbf{x}_0\neq0$,
then the divergence equals 2 regardless of the value of $\mathbf{x}_0$.  






To address this issue, we generalize the divergence of the quadratic score by computing the $L^{2}$ distance between smoothed PDFs. Specifically, we consider
the convolutions between the PDFs and a piecewise smooth function $w$, i.e. $f_{\mathbf{X}}*w$ and $f_{\mathbf{Y}}*w$. The $L^{2}$ divergence between these is given by
\begin{equation}
\intop_{\mathbb{R}^{d}}\left(\left(f_{\mathbf{X}}*w\right)(\mathbf{z})-\left(f_{\mathbf{Y}}*w\right)(\mathbf{z})\right)^{2}\,h({\bf z})\,d{\bf z},\label{eq:L2divergence}
\end{equation}
where $h$ is a weight function needed to ensure the convergence of the integral, which is not necessarily finite with respect to the Lebesgue measure. In some cases, these convolutions have intuitive interpretations. When $w$ is the PDF of a continuous
random variable $\mathbf{W}$, such as a Gaussian random variable, the
convolutions $f_{\mathbf{X}}*w$ and $f_{\mathbf{Y}}*w$ represent
the PDFs of the random variables $\mathbf{X}+\mathbf{W}$ and $\mathbf{Y}+\mathbf{W}$,
respectively; when $w$ is the Heaviside function, i.e. $w(\mathbf{z})=\prod_{j=1}^{d}\I\{z_{j}\geq0\}=:\I\{\mathbf{z}\geq\mathbf{0}\}$,
the convolutions lead to the CDFs.

Recall the example with $f_{\mathbf{Y}}{(\mathbf{z})}=\delta_{\mathbf{0}}(\mathbf{z})$, $f_{\mathbf{X}}(\mathbf{z})=\delta_{\mathbf{x}_{0}}(\mathbf{z})$ and $\mathbf{x}_{0}\neq0$ that we discussed earlier in this section. The weakness that we highlighted with the quadratic score can be addressed by setting $w(\mathbf{z})$ as the Heaviside function in \eqref{eq:L2divergence}. It can be shown that our $L^2$ divergence in \eqref{eq:L2divergence} is then equal to $|\mathbf{x}_{0}|$, which does give credit to more accurate estimates. In addition, the convolutions $f_{\mathbf{X}}*w$ and $f_{\mathbf{Y}}*w$ in \eqref{eq:L2divergence} are piecewise smooth for any
distribution if $w$ is piecewise smooth. This is both theoretically and numerically appealing. 




We can see that \eqref{eq:L2divergence} is nonnegative
and is equal to 0 if $f_{\mathbf{X}}=f_{\mathbf{Y}}$ $\mathcal{L}^{d}$-a.e..
Hence, scores with this divergence function are by definition proper.
Moreover, \eqref{eq:L2divergence} is equal to 0 if and only if
$\left(f_{\mathbf{X}}*w\right)(\mathbf{z})h({\bf z})=\left(f_{\mathbf{Y}}*w\right)(\mathbf{z})h({\bf z})$
$\mathcal{L}^{d}$-a.e.. In this case, if $h$ is nonzero $\mathcal{L}^{d}$-a.e
and $f_{\mathbf{Y}}*w$ uniquely determines any $P_{\mathbf{Y}}\in\mathcal{V}^{d}$ (i.e., $f_{\mathbf{X}}*w=f_{\mathbf{Y}}*w$ $\mathcal{L}^{d}$-a.e.
implies $f_{\mathbf{X}}=f_{\mathbf{Y}}$ $\mathcal{L}^{d}$-a.e.),
then the scores associated with \eqref{eq:L2divergence} are strictly
proper. 

Next, we introduce our $L^{2}$ scoring rule whose divergence
is precisely \eqref{eq:L2divergence}.  The idea
is to view a realization $\mathbf{y}$ as a Dirac delta measure located
at $\mathbf{y}$. The convolution between this Dirac delta measure
and $w$ leads to $w(\mathbf{z}-\mathbf{y})$. Thus, our $L^{2}$ scoring rule can be written as 
\begin{align}
S(P_{{\bf X}},{\bf y};w,h) & =\intop_{\mathbb{R}^{d}}\left(\left(f_{\mathbf{X}}*w\right)(\mathbf{z})-w(\mathbf{z}-\mathbf{y})\right)^{2}h({\bf z})\,d{\bf z}.\label{eq: L2 Score physical}
\end{align}
By expanding the integrand in \eqref{eq: L2 Score physical}, we can
observe that $\intop_{\mathbb{R}^{d}}w{}^{2}(\mathbf{z}-\mathbf{y})h({\bf z})\,d{\bf z}$ only depends on $\mathbf{y}$, hence has the same value regardless
of $P_{\mathbf{X}}$. Therefore, omitting this term has no impact
on the $L^{2}$ divergence in \eqref{eq:L2divergence}. Hence, we can define another scoring rule that is equivalent to \eqref{eq: L2 Score physical} as
\begin{align}
{\color{teal}}S^{'}(P_{\mathbf{X}},\mathbf{y};w,h)= & \intop_{\mathbb{R}^{d}}\left(f_{\mathbf{X}}*w\right)^{2}(\mathbf{z})  h({\bf z})\,d{\bf z} - 2\intop_{\mathbb{R}^{d}}\left(f_{\mathbf{X}}*w\right)(\mathbf{z})w(\mathbf{z}-\mathbf{y})h({\bf z})\,d{\bf z}.\label{eq: L2 score equivalent renormalization physical}
\end{align}
We summarize the results in the following Theorem. The proof is postponed
to the appendix.}

\begin{theorem}\label{thm:main} {\color{\mycolor} Under the conditions
in Assumption \ref{assumption: one-1} in the appendix, the $L^{2}$ scoring rules
in \eqref{eq: L2 Score physical} and \eqref{eq: L2 score equivalent renormalization physical}
are proper. They are strictly proper if $h$ is nonzero $\mathcal{L}^{d}$-a.e. and $f_{\mathbf{y}}*w$ uniquely
characterizes any $P_{\mathbf{Y}}\in\mathcal{V}^{d}$.
Both of the $L^{2}$ scoring rules in (\ref{eq: L2 Score physical})
and (\ref{eq: L2 score equivalent renormalization physical}) have
the divergence function given in \eqref{eq:L2divergence}.}

\end{theorem}
{\color{\mycolor} Recall that our motivation was to use piecewise
smooth functions $w$ to smooth the PDFs via a convolution to improve
upon the quadratic score. From a purely theoretical perspective, the requirement that $w$ is piecewise smooth can be relaxed. For any local Borel measure $w$ and weight
function $h$, if the convolutions are well-defined and the integrals
in \eqref{eq:L2divergence} and \eqref{eq: L2 Score physical} are
finite, then (\ref{eq: L2 Score physical}) is a proper scoring rule.
Similarly, if the convolutions are well-defined and the integrals
in \eqref{eq:L2divergence} and (\ref{eq: L2 score equivalent renormalization physical})
are finite, then (\ref{eq: L2 score equivalent renormalization physical})
is a proper scoring rule. We summarize the technical conditions in
Assumption~\ref{assumption: one-1}. The weight function $h$ can
be anything, but the most convenient form is a PDF, because
PDFs decay rapidly, and so usually guarantee the convergence of (\ref{eq:L2divergence}).
In view of this, we set $h$ to be a PDF in our empirical studies.
The function $h$ can also be used to emphasize regions of interest.
\citet{gneiting2011comparing} put more weight on the tails for inflation
data, and \citet{grushka2017quantile} consider non-constant weight
functions as a way to align the score with the cost function in a
given business context. Moreover, as we show in Section~\ref{subsec:Score Numerical Method},
choosing $h$ to be a PDF is also convenient for numerical computation
of the scores.}

{\color{\mycolor} The proposed $L^{2}$ score is not a single score
but a class of scores. In practice, we want to consider specific $w$
and $h$ so that $\mathcal{V}^{d}$ can be sufficiently general that
important distributions are not excluded. For this purpose, we reiterate
our recommendation that $w$ is a piecewise smooth function and $h$
is a PDF. With these choices, the $L^{2}$ scores can be used with
$\mathcal{V}^{d}$ that covers almost all standard distributions,
e.g. the distributions with bounded PDFs.}

In the literature, there are studies regarding \textit{local scoring
rules, }which depend solely on the PDF or the derivatives of the PDF
purely at the realization $\mathbf{y}$, where a well-known example
is the log score $\log(f_{\mathbf{X}}(\mathbf{y}))$ (see, for example,
\citealt{ehm2012local,parry2012proper}). Because the $L^{2}$ scoring
rules involve integration over $\R^{d}$, they are fundamentally different
from local scoring rules. Hence, the proposed $L^{2}$ scoring rules
do not cover the local scoring rules, and are not exhaustive. The
non-locality of the $L^{2}$ scoring rules is important for our derivation
of the level set scores in Section~\ref{sec:Scoring Functions for Multivariate Level Sets},
because level sets are not local statistical objects in the sense
that their specifications require information regarding the entire
domain. Thus, local scoring functions are not within the scope of
this study, and so we do not discuss them further.

\begin{remark}\label{Remark: Fourier side L2 score}Using the Plancherel
identity, we can also construct $L^{2}$ scoring rules based on characteristic
functions of distributions. In particular, this approach can lead
to the well-known energy score. We discuss this further in the appendix.\end{remark}

\subsection{Particular Cases of $L^{2}$ Scoring Rules\label{subsec: Examples of Scoring Functions}}

In this section, we first demonstrate the generality of our proposed
$L^{2}$ scoring rules by showing that the quadratic score, CRPS,
and MCRPS emerge as special cases of the general framework laid down
in Section~$\ref{subsec: scoring fn}$. We then show how our framework
naturally generates other scoring rules by developing a score based
on lower partial moments. We emphasize that, throughout the examples
in this section, the key issue is the specification {\color{\mycolor}of} $w$.

\subsubsection{Quadratic Score\label{subsec: Density Score}}
\textcolor{\mycolor}{ Let us first note that the convolution between any PDF $f_{\mathbf{X}}(\mathbf{z})$ and   $\delta_{\mathbf{0}}(\mathbf{z})$ is simply $f_{\mathbf{X}}(\mathbf{z})$ itself. In view of this, setting  $w(\mathbf{z})=\delta_{\mathbf{0}}(\mathbf{z})$ in \eqref{eq: L2 score equivalent renormalization physical} gives the following score,}
\begin{equation}
\text{DQS}^{'}(P_{\mathbf{X}},\mathbf{y};h)={\color{\mycolor} S^{'}(P_{\mathbf{X}},\mathbf{y};\delta_{\mathbf{0}},h)}=\intop_{\mathbb{R}^{d}}f_{\mathbf{X}}^{2}(\mathbf{z})h({\bf z})\,d{\bf z}-2f_{\mathbf{X}}(\mathbf{y})h(\mathbf{y}). \label{eq: quadratic score general}
\end{equation}
\textcolor{\mycolor}{When $h\equiv1$, the score in \eqref{eq: quadratic score general} reduces to the well-known quadratic
score.}  For simplicity, from now on, we refer to (\ref{eq: quadratic score general})
as the quadratic score. Note that, for this case, we have used (\ref{eq: L2 score equivalent renormalization physical})
rather than (\ref{eq: L2 Score physical}) because the latter involves
the multiplication of two Dirac delta masses and hence is not well-defined. {\color{\mycolor} The quadratic score can be used with discrete or continuous distributions that have finite PDFs.}


\subsubsection{CRPS \& MCRPS\label{subsec: MCRPS}}

\textcolor{\mycolor} {Recall that $u(\mathbf{z})=\prod_{j=1}^{d}\I\{z_{j}\geq0\}=:\I\{\mathbf{z}\geq\mathbf{0}\}$ is the Heaviside function. As discussed in Section~\ref{subsec: scoring fn}, the convolution of a PDF
$f_{\mathbf{X}}$ and $u$ gives the CDF $F_{\mathbf{X}}$. In fact, we can write }
\begin{align*}
(f_{\mathbf{X}}\ast u)(\mathbf{z}) & =\int_{-\infty}^{\infty}\cdots\int_{-\infty}^{\infty}f_{\mathbf{X}}(s_{1},\ldots,s_{d})\prod_{j=1}^{d}\I\{z_{j}\geq s_{j}\}\,ds_{1}\ldots ds_{d}=:F_{\mathbf{X}}(z_{1},\ldots,z_{d}).
\end{align*}

By Theorem~\ref{thm:main}, {\color{\mycolor} and if we use $w(\mathbf{z}) = u(\mathbf{z})$ in \eqref{eq: L2 Score physical}, we obtain the following proper score:}
\begin{equation}
\text{MCRPS}(P_{\mathbf{X}},\mathbf{y};h) = {\color{\mycolor}S(P_{\mathbf{X}},\mathbf{y};u,h)}=\intop_{\mathbb{R}^{d}}\left(F_{\mathbf{X}}(\mathbf{z})-\I\left\{ \mathbf{z}\geq\mathbf{y}\right\} \right)^{2} h({\bf z})\,d{\bf z}. \label{MCRPS, ex A}
\end{equation}
This score is precisely the MCRPS considered by \citet{gneiting2007strictly}.
For univariate distributions, (\ref{MCRPS, ex A}) is simply the well-known
CRPS for $h\equiv1$, and the threshold weighted CRPS for a general
function $h$ (\citealt{gneiting2011comparing,grushka2017quantile}).

Following (\ref{eq: L2 score equivalent renormalization physical}),
we obtain an equivalent expression for the MCRPS, given by
\begin{equation}
\text{MCRPS}^{'}(P_{\mathbf{X}},\mathbf{y};h)=\intop_{\mathbb{R}^{d}}F_{\mathbf{X}}^{2}(\mathbf{z}) h({\bf z})\,d{\bf z}-2\intop_{\mathbb{R}^{d}}F_{\mathbf{X}}(\mathbf{z})\I\left\{ \mathbf{z}\geq\mathbf{y}\right\} h({\bf z})\,d{\bf z}.\label{eq: MCRPS, renormalization}
\end{equation}
This expression will be useful in our consideration of level sets
in Section~\ref{sec:Scoring Functions for Multivariate Level Sets}.
Both MCRPS and $\text{MCRPS}^{'}$ scores are strictly proper {\color{\mycolor}if $h$ is nonzero $\mathcal{L}^{d}$-a.e.,} because
the CDF uniquely characterizes a distribution. {\color{\mycolor} If $h$ is a PDF, both MCRPS and MCRPS$^{'}$ scores can be used with all distributions.}

\subsubsection{A New Score Based on Lower Partial Moments\label{subsec:Lower-Partial-Moments scoring function}}

\textcolor{\mycolor}{Let $u^{{*k}}(\mathbf{z})\equiv\underbrace{u*u*\cdots*u}_{k\,\,times}(\mathbf{z})=\prod_{j=1}^{d}\frac{1}{k!}z_{j}^{k}\I\{z_{j}\geq0\}$ denote the $k^{\text{th}}$ convolution power of the Heaviside
function, where $k=1,2,3...$. The convolution $f_{\mathbf{X}}*u^{{*k}}$ gives the LPM
of order\textit{ $k$,} which is given by}
\begin{align}
f_{\mathbf{X}}*u^{*k}(\mathbf{z})=\int_{-\infty}^{\infty}\cdots\int_{-\infty}^{\infty}f_{\mathbf{X}}(s_{1},\ldots,s_{d})\prod_{j=1}^{d}\frac{1}{k!}(z_{j}-s_{j})^{k}\I\{z_{j}\geq s_{j}\}\,ds_{1}\ldots ds_{d} := {\rm LPM}_{{\bf X},k}(\mathbf{z}).\label{eq: lower partial moments}
\end{align}
The LPM for univariate distributions has been widely considered for
systemic risk (\citealt{price1982variance}), asset pricing (\citealt{anthonisz2012asset}),
and portfolio management (\citealt{briec2010portfolio}). However,
to the best of our knowledge, the LPM for multivariate distributions
has not been considered in the literature. 

By Theorem 1, {\color{\mycolor}and if we use $w(\mathbf{z}) = u^{*k}(\mathbf{z})$ in \eqref{eq: L2 Score physical}, we obtain} a new proper scoring rule based on the LPM of order~$k$: 
\begin{align*}
\text{LPMS}(P_{\mathbf{X}},\mathbf{y};k,h)=S(P_{\mathbf{X}},\mathbf{y};u^{{*k}},h) & =\intop_{\mathbb{R}^{d}}\left({\rm LPM}_{\mathbf{X},k}(\mathbf{z})-\prod_{j=1}^{d}\frac{1}{k!}(z_{j}-y_{j})^{k}\I\{z_{j}\geq y_{j}\}\right)^{2}h(\mathbf{z})\,d\mathbf{z}.
\end{align*}
We may infer from (\ref{eq: L2 score equivalent renormalization physical}),
the following equivalent score: 
\begin{align}
\text{LPMS}^{'}(P_{\mathbf{X}},\mathbf{y};k,h) & =\intop_{\mathbb{R}^{d}}{\rm LPM}_{\mathbf{X},k}^{2}(\mathbf{z})h({\bf z})\,d\mathbf{z}-2\intop_{\mathbb{R}^{d}}{\rm LPM}_{\mathbf{X},k}(\mathbf{z})\prod_{j=1}^{d}\frac{1}{k!}(z_{j}-y_{j})^{k}\I\{z_{j}\geq y_{j}\}h(\mathbf{z})\,d\mathbf{z}.\label{eq: lower power moments score}
\end{align}
This form of the score will be useful in our consideration of level
sets in Section~\ref{sec:Scoring Functions for Multivariate Level Sets}.
Both $\text{LPMS}$ and $\text{LPMS}^{'}$ scores are strictly proper {\color{\mycolor}
if $h$ is nonzero $\mathcal{L}^{d}$-a.e.,} because the lower partial
moment uniquely characterizes a distribution. A proof is presented
in the appendix. {\color{\mycolor} If $h$ is a PDF with bounded support or exponential
decay, both LPMS and LPMS$^{'}$ scores can be used with distributions that have finite PDFs.}

\section{New Scoring Functions for Level Sets\label{sec:Scoring Functions for Multivariate Level Sets}}

Sometimes we are more interested in a specific region of a distribution
than the entire domain. For example, tails of a distribution are often
of great interest to various applications, with quantiles widely studied
as an important risk measure (see, for example, \citealt{jose2009evaluating,grushka2017quantile}).
Quantiles have also been considered in the context of model averaging
(\citealt{lichtendahl2013better}).

In this paper we consider level sets, which can be viewed as multivariate
generalizations of quantiles (\citealt{abbas2010assessing,cousin2013multivariate}). {\color{\mycolor}Let $g: \mathbb{R}^{d} \rightarrow \mathbb{R}$ be} a univariate function, and $\alpha\in\mathbb{R}$. We define the $\alpha$ level set of $g$ as
\begin{equation}
L(g;\alpha):=\{\mathbf{z}\in\mathbb{R}^{d}:g(\mathbf{z})\geq\alpha\}.\label{eq: general level set}
\end{equation}
For simplicity, we denote $\{\mathbf{z}\in\mathbb{R}^{d}:g(\mathbf{z})\geq\alpha\}$
as $\left\{ g\geq\alpha\right\} $ in the rest of the paper. In the
literature, sometimes $L({\color{teal}{\color{teal}{\color{black}g}}};\alpha)$
is referred to as the upper level set, and $\partial\{L({\color{teal}{\color{teal}{\color{black}g}}};\alpha)\}$
as the level set (see, for example, \citealt{chen2017density}). Our
terminology is consistent with \citet{cadre2006kernel} and \citet{singh2009adaptive}.

More specifically, in this section, we consider $g={\color{teal}{\color{black}f_{\mathbf{Y}}*w}}$,
as {\color{\mycolor}explained} in Section~\ref{sec:A-Framework-for scoring functions},
and {\color{\mycolor}we} derive consistent scoring functions for the level set $L(f_{\mathbf{Y}}*w;\alpha):=\big\{ f_{\mathbf{Y}}*w\geq\alpha\big\}$. As a special case, for a univariate random variable $Y$, {\color{\mycolor}and if  $w$ is the Heaviside function}, the $\alpha$ level set $L(F_{Y};\alpha)$ can be identified with
the $\alpha$ quantile as follows: $q_{Y}(\alpha)\longleftrightarrow[q_{Y}(\alpha),\infty)=L(F_{Y};\alpha).$

Similar to the spirit of a proper scoring rule, we can define scoring
functions for level sets. {\color{\mycolor}Let $l(A,{\bf y};w,h)$ be an $\R$-valued
function  whose arguments consist of a Borel set
$A\subseteq\mathbb{R}^{d}$, an estimate of $L(f_{\mathbf{Y}}*w;\alpha)$,
and a realization ${\bf y}\in\R^{d}$.} We say that $l$ is a consistent
scoring function of the level sets if $\mathbb{E}_{P_{{\bf Y}}}[l(A,\bullet;w)]$
is minimized by $L(f_{\mathbf{Y}}*w;\alpha)$, and strictly consistent
if $L(f_{\mathbf{Y}}*w;\alpha)$ is the only minimizer. In this paper,
to distinguish scores for distributions and level sets, we adopt the
convention that the scores for level sets are termed (consistent)
scoring functions (see, for example, \citealt{fissler2019forecast}).
To the best of our knowledge, scoring functions for general level
sets have not been studied, with the only exceptions being a score
for density level sets $L(f_{{\bf Y}};\alpha)$ and scores for quantiles
of univariate distributions.

In this section, we propose a systematic approach for constructing
scoring functions for level sets. Our inspiration stems from the well-known
result that the CRPS can be decomposed as the integral of the quantile
scores (see, for example, \citealt{gneiting2007strictly,grushka2017quantile}):
\begin{align}
\text{CRPS}(P_{X},y) & :=\intop_{-\infty}^{\infty}\Big(F_{X}(z)-\I\left\{ z\geq y\right\} \Big)^{2}\,dz=2\intop_{0}^{1}\Big(\alpha-\mathbbm{1}\left\{ y<q_{X}(\alpha)\right\} \Big)\Big(y-q_{X}(\alpha)\Big)\,d\alpha,\label{eq: CRPS decomposition}
\end{align}
where the integrand in the last integral is the well-known quantile score, also known as {\color{\mycolor}the ``pinball loss'' or} ``tick-loss''. The equivalence can be established via a change of
variables (see, for example, \citealt{laio2007verification}).
In spite of its simplicity, this algebraic manipulation does not extend
in a straightforward way to other types of scoring functions, and
so, in the multivariate case, it is not trivial to obtain consistent
scoring functions for level sets from proper scoring rules. 

Our idea is to use the ``layer cake representation'' to decompose the scoring
rule $S^{'}(P_{\mathbf{X}},\mathbf{y};w, h)$ in (\ref{eq: L2 score equivalent renormalization physical}) into an integral of scoring functions for level sets of $f_{\mathbf{Y}}*w$.
This provides a unified approach for constructing scoring functions
for different types of level sets, including the scores for density
level sets and quantiles. In addition, our approach provides insight
into the relationship between distributions and their level sets. 

In {\color{\mycolor}the next section}, we formally
describe our approach for constructing scoring functions for level
sets. {\color{\mycolor}We then} consider specific examples {\color{\mycolor}in Section~\ref{subsec: Examples of Scores for Level Sets}}.

\subsection{Scoring Functions for Level Sets\label{subsec: Layer Cake Representation}}

{\color{\mycolor} The key tool for our further developments is the
elementary but useful \textit{layer cake representation}.
Its proof can be found in \citet{lieb2001analysis}, and the name
``layer cake'' refers to the level set structure. Roughly speaking,
the layer cake representation decomposes the $L^{2}$ scoring rules
in (\ref{eq: L2 score equivalent renormalization physical}) into
an integral of ``layers'', which are precisely the level set scores.
Theorem \ref{thm:scoring function for level set} summarizes the theoretical
results. The derivation and proof are postponed to the appendix.}
\begin{theorem}\label{thm:scoring function for level set} Let $\mathcal{V}^{d}$, $h$, and $w$ satisfy
{\color{\mycolor} Assumption~\ref{assumption: one-1} in the appendix}, and $w$ is nonnegative. Let $\alpha>0$, and let $A\subseteq\R^{d}$
be an arbitrary Borel set. Suppose $\left(S'\right)^{\Gamma}\left(A,\bullet;w,h,\alpha\right)$
is finite $\mathcal{L}^{d}$-a.e. and $\left(S'\right)^{\Gamma}\left(A,\bullet;w,h,\alpha\right)f_{{\bf Y}}(\bullet)$
is integrable. {\color{\mycolor} Then the following expression defines a
consistent scoring function for the $\alpha$ level set $L(f_{{\bf Y}}\ast w;\alpha)$:
\begin{align}
(S^{'})^{\Gamma}\Big(A,\mathbf{y};w,h,\alpha\Big) & =\intop_{\mathbb{R}^{d}}\left(\alpha-w(\mathbf{z}-\mathbf{y})\right)\mathbbm{1}\left\{ \mathbf{z}\in A\right\} h({\bf z})\,d{\bf z}.\label{eq:scoring functions for level sets}
\end{align}
If $h$ is nonzero $\mathcal{L}^{d}$-a.e. and} $L(f_{{\bf Y}}\ast w;\alpha)$
equals the closure of $\{f_{{\bf Y}}\ast w>\alpha\}$, it is strictly
consistent. {\color{\mycolor} Moreover, with $A=L(f_{{\bf X}}\ast w;\alpha)$, the integral of (\ref{eq:scoring functions for level sets}) with respect to $\alpha$  is precisely half of the $L^{2}$
score in (\ref{eq: L2 score equivalent renormalization physical}),
i.e., $\frac{1}{2}S^{'}(P_{\mathbf{X}},\mathbf{y};w,h)=\intop_{0}^{\infty}(S^{'})^{\Gamma}\Big(L(f_{{\bf X}}\ast w;\alpha),\mathbf{y};w,h,\alpha\Big)\,d\alpha.$}

\end{theorem}
{\color{\mycolor} The superscript~$^{\Gamma}$ emphasizes that $(S')^{\Gamma}$
is a scoring function for a level set, which is derived from the $L^{2}$
scoring rule $S^{'}(P_{\mathbf{X}},\mathbf{y};w,h)$ in (\ref{eq: L2 score equivalent renormalization physical}).
Theorem \ref{thm:scoring function for level set} provides important
insight. It generalizes the decomposition of the CRPS in (\ref{eq: CRPS decomposition})
into multivariate settings, namely, the $L^{2}$ scoring rule for
$P_{\mathbf{Y}}$ is equal to two times the integral with respect to $\alpha$ of the scoring functions for the level
sets $L(f_{{\bf Y}}\ast w;\alpha)$. Since level sets represent different regions, the level set scores effectively enable us to focus on a specific region by considering specific values of $\alpha$.}

\subsection{Particular Cases of Scoring Functions for Level Sets\label{subsec: Examples of Scores for Level Sets}}

In this section, we first show how our approach allows us to derive
the scoring functions for density level sets (\citealt{hartigan1987estimation,chen2017density}).
We then develop two new scores: the scoring functions for CDF level
sets $L(F_{\mathbf{Y}};\alpha)$ and LPM level sets $L({\rm LPM}_{{\bf Y},k};\alpha)$.
In particular, we show that, in the univariate context, the scoring
functions for CDF level sets induce the full class of quantile scores
studied by \citet{gneiting2011making} and \citet{komunjer2005quasi}.

\subsubsection{Density Level Set Scores \label{subsec: Scoring Function for Density Contour} }

\begin{figure}[H]
\centering{}\caption{\linespread{0.50}\selectfont{}2-D projections of the density level
\textcolor{\mycolor}{sets shown} in Figure~\ref{fig:Illustration level sets 3D}
(a). The numerical values within the plot indicate the value of $\alpha$
for each level set. \textcolor{\mycolor}{Darker coloring corresponds to a higher value of $\alpha$.} \label{fig:Density Level Sets}}
\vspace{-0.2cm}
 \includegraphics[width=0.3\textwidth]{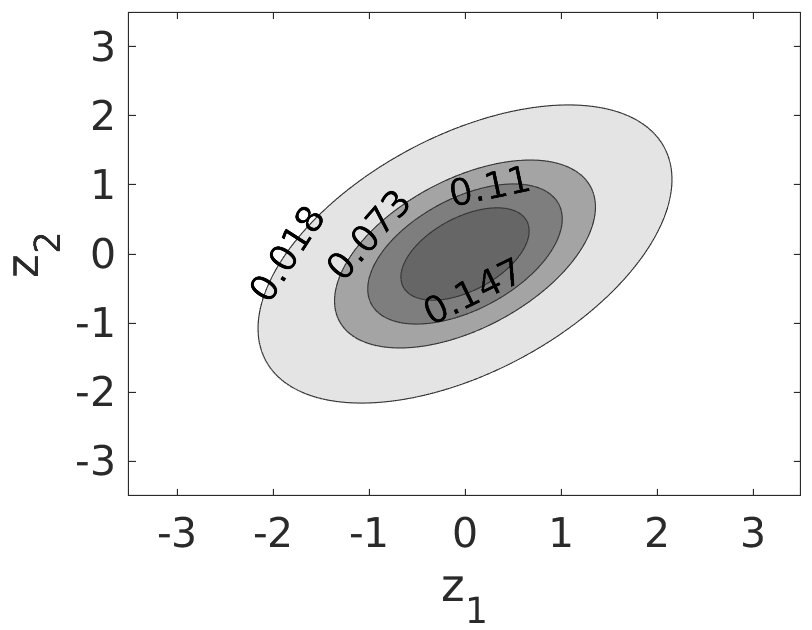}\vspace{-0.8cm}
\end{figure}

Recall that we consider $w(\mathbf{z})\equiv\delta_{\mathbf{0}}(\mathbf{z})$ for the quadratic score $DQS^{'}$ in Section~\ref{subsec: Density Score}. For this $w$, $f_{\mathbf{Y}}*w$ reduces to the PDF $f_{\mathbf{Y}}$, and
by (\ref{eq: general level set}), the $\alpha$ density level set
is defined as $L(f_{\mathbf{Y}};\alpha)=\big\{ f_{\mathbf{Y}}\geq\alpha\big\}$. Figure~\ref{fig:Density Level Sets} shows the 2-D projections of the density level sets of the DGP that were shown in 3-D in Figure \ref{fig:Illustration level sets 3D}(a). {\color{\mycolor} The numerical values within the plot indicate the value $\alpha$ of the PDF $f_{\mathbf{Y}}$ at points on the boundary of the $\alpha$ density level set}.


\begin{figure}[H]
\centering{}\centering{}\caption{\linespread{0.50}\selectfont{}Graphical illustration of the density
level set score in (\ref{eq: scoring fn, density contour, prelim}).
The shaded areas represent the $\alpha$ predictive density level set and $\mathbf{y}$
represents a realization. In (a), $\mathbf{y}$ is inside the level
set, so the density level set score is the Borel measure (induced
by $h$) for the shaded region multiplied by $\alpha$, minus $h(\mathbf{y})$;
in (b), $\mathbf{y}$ is outside the level set, so $\mathbf{y}$ has
no contribution to the score, i.e., the density level set score is
simply the Borel measure for the shaded region multiplied by $\alpha$.{\small{}\label{fig:PDF level set score}}}
 \includegraphics[width=0.6\textwidth]{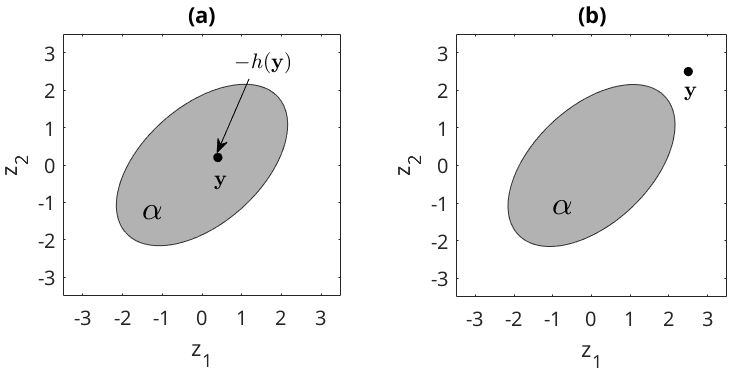}\vspace{-0.8cm}
\end{figure}

\textcolor{\mycolor}{By Theorem~\ref{thm:scoring function for level set}, and using $w(\mathbf{z})\equiv\delta_{\mathbf{0}}(\mathbf{z})$ in \eqref{eq:scoring functions for level sets}}, we obtain the following consistent scoring
\textcolor{\mycolor}{function} for the $\alpha$ density level set,
\begin{align}
{\rm (DQS^{'})^{\Gamma}}\left(A,\mathbf{y};h,\alpha\right) & =\intop_{\mathbb{R}^{d}}\left(\alpha-\delta_{{\bf 0}}(\mathbf{z}-\mathbf{y})\right)\mathbbm{1}\left\{ \mathbf{z}\in A\right\} h({\bf z})\,d{\bf z}\nonumber \\
 & =\alpha\intop_{\mathbb{R}^{d}}\I\big\{\mathbf{z}\in A\big\} h({\bf z})d\mathbf{z}-\I\left\{ \mathbf{y}\in A\right\} h({\bf y}).\label{eq: scoring fn, density contour, prelim}
\end{align}
This is precisely the scoring \textcolor{\mycolor}{function for
density level sets}, also known
as the \textit{excess mass }\textit{\textcolor{\mycolor}{score}} (see, for example, \citealt{hartigan1987estimation,chen2017density}).
The score (\ref{eq: scoring fn, density contour, prelim}) has an
intuitive graphical interpretation, which we illustrate in Figure~\ref{fig:PDF level set score}.

\begin{figure}[H]
\centering{}\caption{\linespread{0.50}\selectfont{}For $\alpha=0.018$, comparison of
the density level set of the DGP in (a) and the following three misspecified
distributions: (b) misspecified means $[-0.4,-0.4]$; (c) misspecified
variances $[0.6,0.6]$; and (d) misspecified covariance 0.5. The dotted
lines in (b), (c) and (d) indicate the density level set of the DGP
in (a). The numerical value at the top right corner of each plot is the average score\textcolor{\mycolor}{\footnotesize{}$(\times10^{6})$} computed
in the simulation study of Section~\ref{subsec:Simulation-Results}.\label{fig:Figure 4 PDF}}
 \includegraphics[width=1\textwidth]{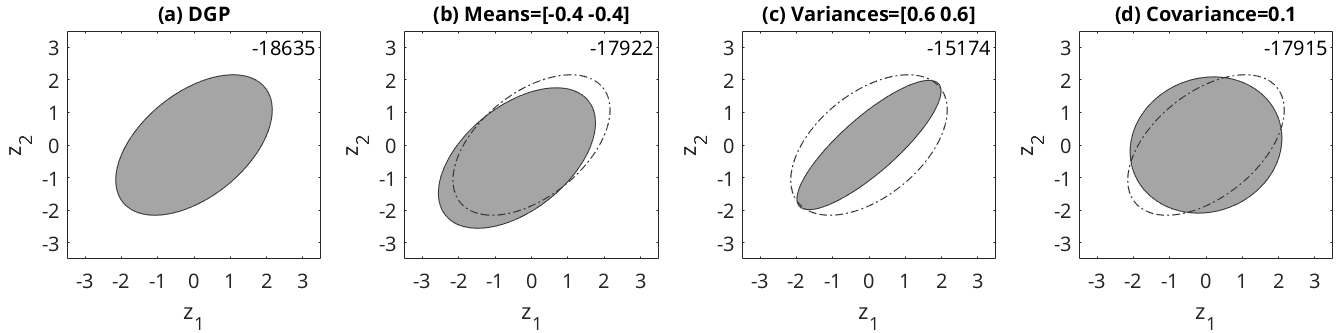}\vspace{-0.6cm}
\end{figure}

We note that competing density forecasts can be compared in terms
of their accuracy for different regions of the density by evaluating
the level sets for different choices of $\alpha$. Focusing on $\alpha=0.018$,
Figure~\ref{fig:Figure 4 PDF} compares the density level set of
the DGP and that of three misspecified distributions, along with the
mean of the density level set scores computed in our simulation study
of Section~\ref{subsec:Simulation-Results}. The DGP receives the
lowest mean score, which shows the consistency of the score. 

\subsubsection{New CDF Level Set Score\label{subsec: Scoring Function for CDF level sets} }

\begin{figure}[H]
\centering{}\caption{\linespread{0.50}\selectfont{}2-D projections of the CDF level \textcolor{\mycolor}{sets
shown} in Figure~\ref{fig:Illustration level sets 3D} (b). The numerical
values within the plot indicate the value of $\alpha$ for each level
set. {\color{\mycolor} Darker coloring corresponds to a higher value of $\alpha$}. {\small{}\label{fig:Density Level Sets-1}}}
\vspace{-0.3cm}
 \includegraphics[width=0.3\textwidth]{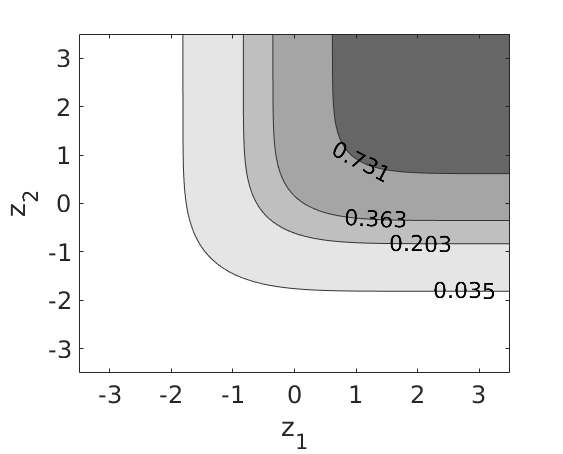}\vspace{-0.8cm}
\end{figure}

Recall that we consider $w(\mathbf{z})=u(\mathbf{z})=\I\{\mathbf{z}\geq\mathbf{0}\}$
in Section~\ref{subsec: MCRPS} to derive $\text{MCRPS}$ and $\text{MCRPS}{}^{'}$. For this $w$, $f_{\mathbf{Y}}*w$ reduces to the CDF $F_{\mathbf{Y}}$, and
by (\ref{eq: general level set}), the $\alpha$ CDF level set is
defined as $L(F_{\mathbf{Y}};\alpha)=\big\{ F_{\mathbf{Y}}\geq\alpha\big\}$. Figure~\ref{fig:Density Level Sets-1} shows the 2-D projections of the CDF level sets of the DGP that were shown in 3D in Figure~\ref{fig:Illustration level sets 3D}(b). {\color{\mycolor} The numerical values within the plot indicate the value $\alpha$ of the CDF $F_{\mathbf{Y}}$ at points on the boundary of the $\alpha$ CDF level set, i.e. $\alpha$ is a probability}.


\begin{figure}[H]
\centering{}\caption{\label{fig:CDF level set score}\linespread{0.50}\selectfont{}Graphical
illustration of the CDF level set score in (\ref{eq: score for isoprobability contour}).
The curves represent the $\alpha$ CDF level set and $\mathbf{y}$
represents a realization, and $\alpha$ or $\alpha-1$ is the weight
for the corresponding shaded region. Three scenarios are presented
depending on the position of $\mathbf{y}$ relative to the level set. The CDF
level set score is the weighted sum of the Borel measures (induced
by $h$) for the shaded regions.}

\includegraphics[width=0.75\textwidth]{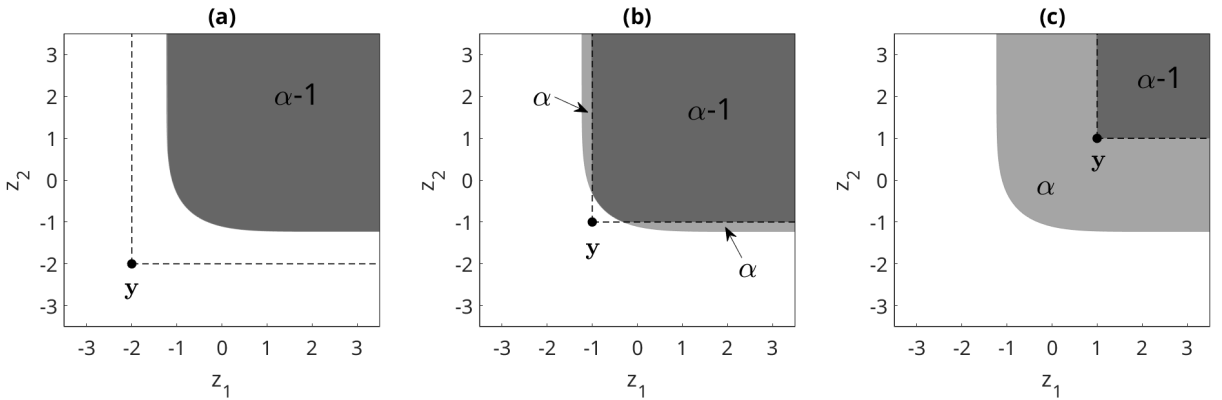}\vspace{-0.6cm}
\end{figure}

\textcolor{\mycolor}{By Theorem~\ref{thm:scoring function for level set}, and using $w(\mathbf{z})=u(\mathbf{z})$
in (\ref{eq:scoring functions for level sets})},
we obtain a consistent scoring function for the $\alpha$ CDF level set,
\begin{equation}
(\text{MCRPS}{}^{'})^{\Gamma}\left(A,\mathbf{y};h,\alpha\right)={\color{teal}{\color{black}\intop_{\mathbb{R}^{d}}\left(\alpha-\I\big\{\mathbf{z}\geq\mathbf{y}\big\}\right)\I\big\{\mathbf{z}\in A\big\} h({\bf z})\,d{\bf z}.}}\label{eq: score for isoprobability contour}
\end{equation}
The score (\ref{eq: score for isoprobability contour}) also has an
intuitive graphical interpretation, which we illustrate in Figure~\ref{fig:CDF level set score}.

\begin{figure}[H]
\centering{}\caption{\linespread{0.50}\selectfont{}For $\alpha=0.035$, comparison of
the CDF level set of the DGP in (a) and the following three misspecified
distributions: (b) misspecified means $[-0.4,-0.4]$; (c) misspecified
variances $[0.6,0.6]$; and (d) misspecified covariance 0.5. The dotted
lines in (b), (c) and (d) indicate the CDF level set of the DGP in
(a). The numerical value at the top right corner of each plot is the average score {\color{\mycolor}\footnotesize{}$(\times10^{5})$} computed
in the simulation study of Section~\ref{subsec:Simulation-Results}.\label{fig:Figure 4 CDF}}
 \includegraphics[width=1\textwidth]{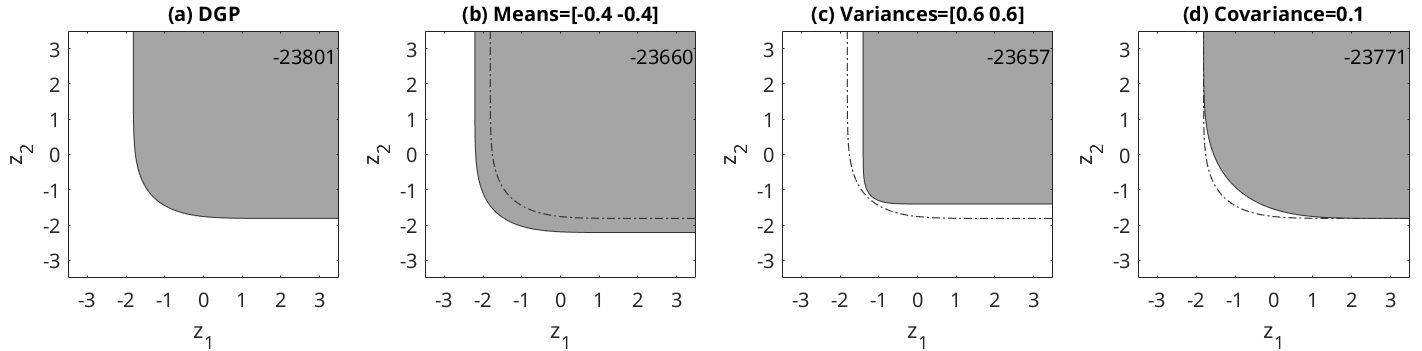}\vspace{-0.6cm}
\end{figure}

Focusing on $\alpha=0.035$, Figure~\ref{fig:Figure 4 CDF} compares
the CDF level set of the DGP and that of three misspecified distributions,
along with the average scores computed in our simulation study of
Section~\ref{subsec:Simulation-Results}. The DGP has the lowest
average score, showing the consistency of the CDF level set score. 

\begin{remark}\label{Remark: quantile score}When $d=1$, (\ref{eq: score for isoprobability contour})
is equivalent to the full class of quantile scores. To see this, let
us first identify a quantile estimate $q(\alpha)$ with a Borel set
$A=[q(\alpha),\infty)$ and then add to (\ref{eq: score for isoprobability contour})
the term $(1-\alpha)\intop_{-\infty}^{\infty}\I\big\{ z\geq y\big\} h(z)\,dz$,
which depends only on $y$, as follows: 
\begin{align*}
 & \intop_{-\infty}^{\infty}\left(\alpha-\I\big\{ z\geq y\big\}\right)\I\big\{ z\geq q(\alpha)\big\} h({\bf z})\,d{\bf z}+(1-\alpha)\intop_{-\infty}^{\infty}\I\big\{ z\geq y\big\} h({\bf z})\,d{\bf z}\\
 & =\alpha\mathbf{\intop_{-\infty}^{\infty}}\left(\I\big\{ z\geq q(\alpha)\big\}-\I\big\{ z\geq y\big\}\right)h({\bf z})\,d{\bf z}+\intop_{-\infty}^{\infty}\left(1-\I\big\{ z\geq q(\alpha)\big\}\right)\I\big\{ z\geq y\big\} h({\bf z})\,d{\bf z}\\
 & =\alpha\intop_{q(\alpha)}^{y}\lambda\left(d{\bf z}\right)-\I\{y<q(\alpha)\}\intop_{q(\alpha)}^{y}h({\bf z})\,d{\bf z}\\
 & =\big(\alpha-\I\{y<q(\alpha)\}\big)\Big(H(y)-H\big(q(\alpha)\big)\Big),
\end{align*}
where $H$ is an anti-derivative of $h$, hence is non-decreasing
(since $h$ is non-negative). The expression in the final line is
the full class of quantile scores (see, for example, \citealt{gneiting2011making,komunjer2005quasi}).\end{remark}

\subsubsection{New LPM Level Set Score}

Recall that we consider $w(\mathbf{z})=u^{*k}(\mathbf{z})=\prod_{j=1}^{d}\frac{1}{k!}z_{j}^{k}\I\{z_{j}\geq0\}$
in Section~\ref{subsec:Lower-Partial-Moments scoring function} to
derive the novel $L^{2}$ scoring rules $\text{LPMS}$ and $\text{LPMS}{}^{'}$.
For this $w$, $f_{\mathbf{Y}}*w$ reduces to the $k$-th lower partial
moment function of the distribution $P_{\mathbf{Y}}$, and by (\ref{eq: general level set})
the $\alpha$ level set of the $k$-th LPM is defined by $L({\rm LPM}_{\mathbf{Y},k};\alpha)=\left\{ {\rm LPM}_{\mathbf{Y},k}\geq\alpha\right\} $.

\begin{figure}[H]
\centering{}\caption{\linespread{0.50}\selectfont{}2-D projections of the LPM level sets
for $k=1$ \textcolor{\mycolor}{shown} in Figure~\ref{fig:Illustration level sets 3D}
(c). The numerical values within the plot indicate the value of $\alpha$
for each level set. {\color{\mycolor}Darker coloring corresponds to a higher value of $\alpha$. {\small{}{}{}{}{}{}\label{fig:LPM Level Sets}}}}
\vspace{-0.3cm}
 \includegraphics[width=0.3\textwidth]{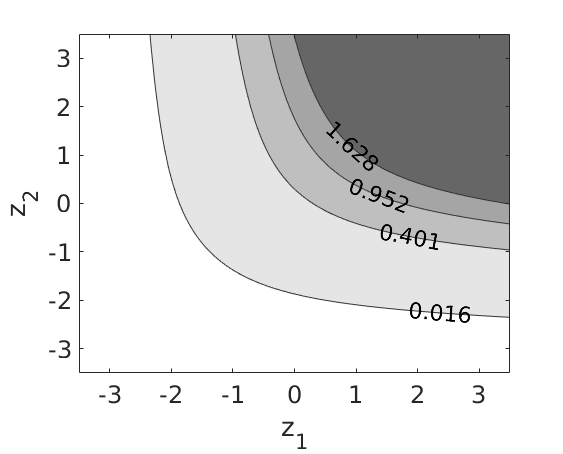}\vspace{-0.8cm}
\end{figure}

{\color{\mycolor}Figure~\ref{fig:LPM Level Sets} shows} the 2-D projections of the LPM level sets for $k=1$ of the DGP that were shown in 3D in Figure~\ref{fig:Illustration level sets 3D}(c). {\color{\mycolor} The numerical values within the plot indicate the value $\alpha$ of the ${\rm LPM}_{\mathbf{Y},1}$ at points on the boundary of the $\alpha$ level set of the ${\rm LPM}_{\mathbf{Y},1}$}.


{\color{\mycolor}By Theorem~\ref{thm:scoring function for level set}, and using $w(\mathbf{z})=u^{*k}(\mathbf{z})$
in (\ref{eq:scoring functions for level sets})},
we obtain a consistent scoring function for the $\alpha$ level set of the
$k$-th LPM: 
\begin{align}
(\text{LPMS}{}^{'})^{\Gamma} & \left(A,\mathbf{y};h,\alpha,k\right):=\intop_{\mathbb{R}^{d}}\left(\alpha-\I\left\{ \mathbf{z}\geq\mathbf{y}\right\} \prod_{j=1}^{d}\frac{(z_{j}-y_{j})^{k}}{k!}h(\mathbf{z})\right)\I\left\{ \bm{{z}}\in A\right\} \,d\mathbf{z}.\label{eq: LPM level set score}
\end{align}

\begin{figure}[H]
\centering{}\caption{\linespread{0.50}\selectfont{}For $k=1$ and $\alpha=0.016$, comparison
of the LPM level set of the DGP in (a) and the following three misspecified
distributions: (b) misspecified means $[-0.4,-0.4]$; (c) misspecified
variances $[0.6,0.6]$; and (d) misspecified covariance 0.5. The dotted
lines in (b), (c) and (d) indicate the LPM level set of the DGP in
(a). The numerical value at the top right corner of each plot is the average score {\color{\mycolor} \footnotesize{}$(\times10^{5})$} in the simulation
study in Section~\ref{subsec:Simulation-Results}.\label{fig: Figure 4 LPM}}
 \includegraphics[width=1\textwidth]{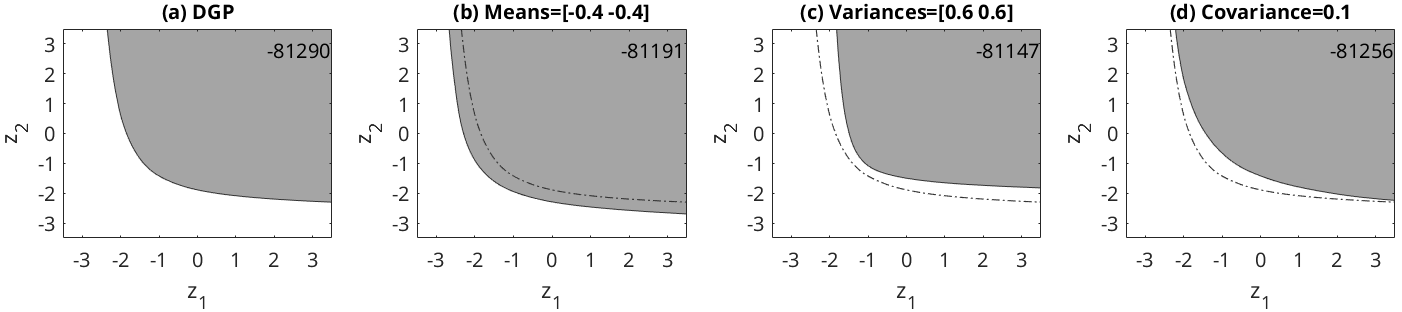}\vspace{-0.8cm}
\end{figure}

\noindent We note that, unlike the density and CDF level set scores,
the LPM level set score does not have an intuitive graphical interpretation.
For $k=1$ and $\alpha=0.016$, Figure~\ref{fig: Figure 4 LPM} compares
the LPM level set of the DGP and that of three misspecified distributions,
along with the average scores computed in the simulation study of
Section~\ref{subsec:Simulation-Results}. The DGP receives the lowest
average score, showing the consistency of the LPM level set score.

\noindent \begin{remark}We remark that our level set scores are related
to the exhaustive scores for the ``set-valued functionals'' studied
by \citet{fissler2019forecast}. They utilize a ``bottom-up'' approach
that constructs the scores from identification functions, while we
consider a ``top-down'' approach, where the level set scores are obtained
through the decomposition of the $L^{2}$ scores. In some cases, similar results
can be obtained from the two approaches. For example, the CDF level
set can be identified with a Vorob'ev quantile of the random set $\left\{ \mathbf{z}\geq\mathbf{Y}\right\} $.
As a result, for this particular example, our CDF level set score
in (\ref{eq: score for isoprobability contour}) coincides with the
Vorob'ev quantile scores in (5.2) in \citet{fissler2019forecast}.\end{remark}

\section{{\color{\mycolor} Numerical Illustration}\label{sec:Simulation-Study}}

In this section, we show how to compute the different scores we have previously introduced using a simulation-based method. {\color{\mycolor} We also perform a simulation study to demonstrate the use of these scores. We implement the $L^{2}$ scoring rules for multivariate distributions, which we presented in Section~\ref{sec:A-Framework-for scoring functions}. These are the ${\rm DQS^{'}}$ in (\ref{eq: quadratic score general}), the ${\rm MCRPS^{'}}$ in (\ref{eq: MCRPS, renormalization}), and the ${\rm LPMS^{'}}$ in (\ref{eq: lower power moments score}). For each $L^{2}$ scoring rule, we also implement the corresponding level set scoring functions described in Section~\ref{sec:Scoring Functions for Multivariate Level Sets}: the $({\rm DQS^{'}})^{\Gamma}$ for density level sets in (\ref{eq: scoring fn, density contour, prelim}), the $({\rm MCRPS^{'}})^{\Gamma}$ for CDF level sets in (\ref{eq: score for isoprobability contour}), and the $({\rm LPMS^{'}})^{\Gamma}$ for LPM level sets in (\ref{eq: LPM level set score}) for $k=1$. Recall that all these $L^{2}$ scoring rules and level set scoring functions are particular cases of our score $S^{'}$ in (\ref{eq: L2 score equivalent renormalization physical}) and $(S')^{\Gamma}$ in (\ref{eq:scoring functions for level sets}), respectively.}


\subsection{Computation of Scores\label{subsec:Score Numerical Method}}

Since the computation of the $L^{2}$ scoring rule $S^{'}$ in (\ref{eq: L2 score equivalent renormalization physical})
and the level set score $(S')^{\Gamma}$ in (\ref{eq:scoring functions for level sets}) essentially involves integrations over $\mathbb{R}^{d}$, one can adopt any numerical integration approach. {\color{\mycolor} If we assume $h$ is a probability measure, the integrals can be simply computed
 as expectations  with respect to the distribution
induced by $h$.} More specifically, let $\left\{ \mathbf{z}^{i}\right\} _{i=1}^{M}$
be a random sample drawn from the distribution
with PDF $h$. The $L^{2}$ scoring rule in (\ref{eq: L2 score equivalent renormalization physical})
can be computed as 
\begin{align}
{\color{teal}}\overline{S^{'}}(P_{\mathbf{X}},\mathbf{y};w,h)=\frac{1}{M} & \sum_{i=1}^{M}\left\{ \left(f_{\mathbf{X}}*w\right)^{2}(\mathbf{z}^{i})-2\left(f_{\mathbf{X}}*w\right)(\mathbf{z}^{i})w(\mathbf{z}^{i}-\mathbf{y})\right\} ,\label{eq: L2 score computation}
\end{align}
and the level set score in (\ref{eq:scoring functions for level sets})
can be computed as\footnotemark{}\footnotetext{\linespread{0.50}\selectfont{} For the quadratic score $\text{DQS}^{'}(P_{\mathbf{X}},\mathbf{y};h)$ in (\ref{eq: quadratic score general}) and the density level set score in (\ref{eq: scoring fn, density contour, prelim}), the terms $f({\bf y})h({\bf y})$ and $\I\left\{\bf y\in A\right\}h({\bf y})$ do not require numerical integration.}
\begin{align}
{\color{teal}}\overline{(S^{'})^{\Gamma}}\left\{ A,\mathbf{y};w,h,\alpha\right\}  & =\frac{1}{M}\sum_{i=1}^{M}\left\{ \left(\alpha-w(\mathbf{z}^{i}-\mathbf{y})\right)\I\left\{ \mathbf{z}^{i}\in A\right\} \right\}. \label{eq: level set score computation}
\end{align}
Following the Central Limit Theorem, the error of this numerical approach
decays at the rate of $O\left(M^{-\frac{1}{2}}\right)$ (see \citealt{caflisch1998monte,robert2013monte}).
It is quick to simulate random \textcolor{\mycolor}{values} using methods
such as MCMC, and the numerical integration boils down to matrix operation,
which can be computed swiftly using parallel computing. In our numerical
studies, we primarily chose $h$ to be a uniform distribution
on the bounded region $[a,b]^{d}$. 

{\color{\mycolor} When comparing different models, it is the relative
difference between the scores that is important, not the absolute
values. In reporting numerical results in this paper, we set one model
as the reference, and then subtract the score of this reference model
from the score for each other model. This leads to a clearer presentation
of the results, and has no impact on the ranking of the performance
of the models.}

\subsection{Simulation Study\label{subsec:Simulation-Results}}

We use the scores to compare the fit of candidate distributions and
their level sets to simulated data. {\color{\mycolor}We generate $T=2\times10^{5}$
observations from the DGP of our running example}, i.e.,
a bivariate normal distribution with zero means, unit variances, and
covariance 0.5.

For CDF level sets, $\alpha$ has range $[0,1]$ for all distributions.
For density and LPM level sets, however, $\alpha$ does not have a
universal range for all distributions. To select a set of values of
$\alpha$ for each type of level set, we first recorded the value
of $\alpha$ for the level set on which each of the $T$ simulated
observations was located. We then chose $\alpha$ to be the $0.1,0.2,...,0.9$
quantiles of these values. The resultant values of $\alpha$ for each
type of level set are presented in the second rows of Tables~\ref{tab:Density level set_simulation}-\ref{tab:CDF level set_simulation Gaussian 2}.

We compared 13 candidate distributions in terms of their fit to the
simulated observations. All were bivariate normal. One had no misspecification,
i.e. it was the DGP; four had the following misspecified means {\color{\mycolor} {[}-2,-2{]},
{[}-0.4,0.4{]}, {[}0.4,0.4{]}, {[}2,2{]}, }and no misspecification
in the variances and covariance; four had the following misspecified
variances {[}0.6,0.6{]}, {[}0.8,0.8{]}, {[}1.2,1.2{]}, {[}1.4,1.4{]},
and no misspecification in the means and covariance; four had the
following misspecified covariance {\color{\mycolor}{}-0.5, 0.1, 0.7,
and 0.9}, and no misspecification in the means and variances. We calculated
the scores using (\ref{eq: L2 score computation}) and (\ref{eq: level set score computation}),
as described in Section~\ref{subsec:Score Numerical Method}, with
$h$ chosen to be the PDF of a uniform distribution on $[-3,3]^{2}$,
and $M=10^{5}$ values sampled from this distribution. For each of
the $T$ observations, we computed the scoring function for each type
of level set and the corresponding $L^{2}$ scoring rule. {\color{\mycolor} We
chose the DGP as the reference, and subtracted its scores from the scores for each
candidate.} For each method, lower values of this score difference
are preferable.

Table~\ref{tab:Density level set_simulation} presents the score
differences for the density level set score ${\rm (DQS^{'})^{\Gamma}}$
and its corresponding $L^{2}$ scoring rule ${\rm DQS^{'}}$; Table~\ref{tab:CDF level set_simulation}
presents the results for the CDF level set score ${\rm (MCRPS^{'})^{\Gamma}}$
and the scoring rule ${\rm MCRPS^{'}}$; and Table~\ref{tab:LPM level set_simulation}
presents the results for the LPM level set score ${\rm (LPMS^{'})^{\Gamma}}$
and the scoring rule ${\rm LPMS^{'}}$. In each table, all entries
for the misspecified candidates are strictly positive, indicating
that the scores are able to identify the DGP. This supports our assertion
that the new scoring rule ${\rm LPMS^{'}}$ is proper and that ${\rm (DQS^{'})^{\Gamma}}$,
${\rm (MCRPS^{'})^{\Gamma}}$, and ${\rm (LPMS^{'})^{\Gamma}}$ are
consistent scoring functions for density level sets, CDF level sets,
and LPM level sets, respectively. It also supports the more general
theoretical result that the $L^{2}$ scoring rule $S^{'}$ in (\ref{eq: L2 score equivalent renormalization physical})
is proper, and $(S^{'})^{\Gamma}$ in (\ref{eq:scoring functions for level sets})
is consistent for level sets.

Next, we show that the properness of our $L^{2}$ scoring rules and the consistency of the level set scores also hold for other choices of
$h$. Tables~\ref{tab:CDF level set_simulation Gaussian} and \ref{tab:CDF level set_simulation Gaussian 2}
present results for the density level set score ${\rm (DQS^{'})^{\Gamma}}$
and its corresponding $L^{2}$ scoring rule ${\rm DQS^{'}}$ with
$h$ chosen as the PDF of the bivariate normal distributions $\mathcal{N}(\mathbf{0},\mathbf{I}_{2})$
and $\mathcal{N}(\mathbf{0},2\times\mathbf{I}_{2})$, respectively, where $\mathbf{I}_{2}$
denotes the identity matrix. To save space, we do not present the
results for the other scores. It can be seen that all values for the
misspecified candidates are strictly positive, indicating that the
DGP receives the lowest scores. {\color{\mycolor} Tables \ref{tab:Density level set_simulation}-\ref{tab:CDF level set_simulation Gaussian 2}
were based on a single run of the simulation study (with $T=2\times10^{5}$
observations). We repeated the study 1000 times,
and found the DGP almost always received the lowest scores. This shows
that the proposed scores can reliably distinguish the DGP from the
misspecified distributions.}

To conclude this section, we briefly discuss the role of $h$ as a
weight function. In the literature, there are the two major weighting
schemes, threshold weighting and quantile weighting, which have been
applied to scores for univariate distributions (\citealt{gneiting2011comparing,grushka2017quantile}).
The function $h$ we consider in our paper enables threshold weighting.
For example, the uniform PDF considered in Tables~\ref{tab:Density level set_simulation}-\ref{tab:LPM level set_simulation}
ignore all the regions outside $[-3,3]^{d}$, while the bivariate
normal PDFs considered in Tables~\ref{tab:CDF level set_simulation Gaussian}-\ref{tab:CDF level set_simulation Gaussian 2}
assign more weight to the center of the distributions. The function
$h$ cannot enable quantile weighting, and so, for example, cannot
be used to emphasize regions of a univariate distribution related
to VaR for a chosen probability level. For this, one could use either
the quantile score or the quantile weighted CRPS. Our level set scores
can be viewed as generalizations of the quantile score to the multivariate
setting, which enable a complementary and useful way to emphasize
regions of interest.
The consistency of the level set score holds for any $h$  {\color{\mycolor}that is a PDF}, which means the regions of interest
in the level set scores are governed by the level sets and not by
$h$. In practice, one could consider different $h$ to compare two
misspecified estimates. If one candidate consistently outperforms
the other for a large family of $h$, then we could conclude that
the former is \textit{dominating} the latter. For example, \citet{ehm2016quantiles}
essentially consider different $h$ to produce a so-called ``Murphy
diagram'' to graphically compare two quantile estimates for univariate
distributions.

\begin{table}[H]
{\footnotesize{}\caption{\linespread{0.50}\selectfont{}{\footnotesize{}For the simulated data,
comparison of candidate distributions using the density level set
score ${\rm (DQS^{'})^{\Gamma}}$ and $L^{2}$ score ${\rm DQS^{'}}$
$(\times10^{6})$. The scores were computed with $h$ chosen as the
PDF of the uniform distribution on $[-3,3]^{2}$. Each value is the
score for a candidate distribution minus the score for DGP. Lower
values are better.}\label{tab:Density level set_simulation}}
\vspace{0.1cm}
}\linespread{0.40}\selectfont{}{\footnotesize{} }%
\begin{tabular*}{1\textwidth}{@{\extracolsep{\fill}}>{\raggedright}p{0.1cm}cccccccccccc}
\toprule 
 &  & \multicolumn{9}{c}{{\footnotesize{}${\rm (DQS^{'})^{\Gamma}}$}} &  & {\footnotesize{}${\rm DQS^{'}}$}\tabularnewline
\midrule 
 & $\alpha$ & {\footnotesize{}0.018} & {\footnotesize{}0.037} & {\footnotesize{}0.055} & {\footnotesize{}0.073} & {\footnotesize{}0.092} & {\footnotesize{}0.110} & {\footnotesize{}0.129} & {\footnotesize{}0.147} & {\footnotesize{}0.165} &  & \tabularnewline\addlinespace[0.1cm]
\midrule 
\multicolumn{2}{c}{{\footnotesize{}DGP}} & {\footnotesize{}0} & {\footnotesize{}0} & {\footnotesize{}0} & {\footnotesize{}0} & {\footnotesize{}0} & {\footnotesize{}0} & {\footnotesize{}0} & {\footnotesize{}0} & {\footnotesize{}0} &  & {\footnotesize{}0}\tabularnewline\addlinespace[0.1cm]
\multicolumn{4}{l}{{\footnotesize{}Misspecified means}} &  &  &  &  &  &  &  &  & \tabularnewline
 & {\footnotesize{}{[}-2,-2{]}} & {\footnotesize{}13297} & {\footnotesize{}13840} & {\footnotesize{}13395} & {\footnotesize{}12558} & {\footnotesize{}11502} & {\footnotesize{}10082} & {\footnotesize{}7611} & {\footnotesize{}5281} & {\footnotesize{}2802} &  & {\footnotesize{}3497}\tabularnewline
 & {\footnotesize{}{[}-0.4,-0.4{]}} & {\footnotesize{}713} & {\footnotesize{}950} & {\footnotesize{}1071} & {\footnotesize{}1055} & {\footnotesize{}1026} & {\footnotesize{}1061} & {\footnotesize{}774} & {\footnotesize{}566} & {\footnotesize{}341} &  & {\footnotesize{}280}\tabularnewline
 & {\footnotesize{}{[}0.4, 0.4{]}} & {\footnotesize{}657} & {\footnotesize{}941} & {\footnotesize{}1096} & {\footnotesize{}1024} & {\footnotesize{}919} & {\footnotesize{}908} & {\footnotesize{}614} & {\footnotesize{}502} & {\footnotesize{}299} &  & {\footnotesize{}261}\tabularnewline
 & {\footnotesize{}{[}2, 2{]}} & {\footnotesize{}13270} & {\footnotesize{}13780} & {\footnotesize{}13326} & {\footnotesize{}12432} & {\footnotesize{}11349} & {\footnotesize{}9966} & {\footnotesize{}7500} & {\footnotesize{}5109} & {\footnotesize{}2602} &  & {\footnotesize{}3459}\tabularnewline\addlinespace[0.1cm]
\multicolumn{5}{l}{{\footnotesize{}Misspecified variances}} &  &  &  &  &  &  &  & \tabularnewline
 & {\footnotesize{}{[}0.6,0.6{]}} & {\footnotesize{}3461} & {\footnotesize{}2342} & {\footnotesize{}1447} & {\footnotesize{}899} & {\footnotesize{}608} & {\footnotesize{}633} & {\footnotesize{}542} & {\footnotesize{}725} & {\footnotesize{}1093} &  & {\footnotesize{}1981}\tabularnewline
 & {\footnotesize{}{[}0.8,0.8{]}} & {\footnotesize{}375} & {\footnotesize{}256} & {\footnotesize{}134} & {\footnotesize{}57} & {\footnotesize{}30} & {\footnotesize{}118} & {\footnotesize{}118} & {\footnotesize{}252} & {\footnotesize{}354} &  & {\footnotesize{}145}\tabularnewline
 & {\footnotesize{}{[}1.2,1.2{]}} & {\footnotesize{}127} & {\footnotesize{}101} & {\footnotesize{}38} & {\footnotesize{}3} & {\footnotesize{}106} & {\footnotesize{}370} & {\footnotesize{}496} & {\footnotesize{}797} & {\footnotesize{}288} &  & {\footnotesize{}86}\tabularnewline
 & {\footnotesize{}{[}1.4,1.4{]}} & {\footnotesize{}429} & {\footnotesize{}250} & {\footnotesize{}107} & {\footnotesize{}213} & {\footnotesize{}566} & {\footnotesize{}1487} & {\footnotesize{}1601} & {\footnotesize{}797} & {\footnotesize{}288} &  & {\footnotesize{}220}\tabularnewline\addlinespace[0.1cm]
\multicolumn{5}{l}{{\footnotesize{}Misspecified covariance}} &  &  &  &  &  &  &  & \tabularnewline
 & {\footnotesize{}-0.5} & {\footnotesize{}3709} & {\footnotesize{}3649} & {\footnotesize{}3204} & {\footnotesize{}2622} & {\footnotesize{}1912} & {\footnotesize{}1342} & {\footnotesize{}641} & {\footnotesize{}336} & {\footnotesize{}172} &  & {\footnotesize{}693}\tabularnewline
 & {\footnotesize{}0.1} & {\footnotesize{}720} & {\footnotesize{}726} & {\footnotesize{}619} & {\footnotesize{}478} & {\footnotesize{}319} & {\footnotesize{}283} & {\footnotesize{}215} & {\footnotesize{}319} & {\footnotesize{}288} &  & {\footnotesize{}154}\tabularnewline
 & {\footnotesize{}0.7} & {\footnotesize{}458} & {\footnotesize{}376} & {\footnotesize{}307} & {\footnotesize{}207} & {\footnotesize{}137} & {\footnotesize{}157} & {\footnotesize{}107} & {\footnotesize{}105} & {\footnotesize{}174} &  & {\footnotesize{}99}\tabularnewline
 & {\footnotesize{}0.9} & {\footnotesize{}3766} & {\footnotesize{}2890} & {\footnotesize{}2090} & {\footnotesize{}1631} & {\footnotesize{}1303} & {\footnotesize{}1239} & {\footnotesize{}1139} & {\footnotesize{}1199} & {\footnotesize{}1324} &  & {\footnotesize{}1334}\tabularnewline
\bottomrule
\end{tabular*}{\footnotesize\par}

{\small{}\vspace{-0.6cm}
}{\small\par}
\end{table}

\begin{table}[H]
{\footnotesize{}\caption{\linespread{0.50}\selectfont{}{\footnotesize{}For the simulated data,
comparison of candidate distributions using the CDF level set score
${\rm (MCRPS^{'})^{\Gamma}}$ and $L^{2}$ score ${\rm MCRPS^{'}}$
$(\times10^{5})$. The scores were computed with $h$ chosen as the
PDF of the uniform distribution on $[-3,3]^{2}$. Each value is the
score for a candidate distribution minus the score for DGP. Lower
values are better. }\label{tab:CDF level set_simulation}}
\vspace{0.1cm}
}\linespread{0.40}\selectfont{}{\footnotesize{} }%
\begin{tabular*}{1\textwidth}{@{\extracolsep{\fill}}>{\raggedright}p{0.1cm}cccccccccccc}
\toprule 
 &  & \multicolumn{9}{c}{{\footnotesize{}${\rm (MCRPS^{'})^{\Gamma}}$}} &  & {\footnotesize{}${\rm MCRPS^{'}}$}\tabularnewline
\midrule 
 & $\alpha$ & {\footnotesize{}0.035} & {\footnotesize{}0.082} & {\footnotesize{}0.138} & {\footnotesize{}0.203} & {\footnotesize{}0.277} & {\footnotesize{}0.363} & {\footnotesize{}0.462} & {\footnotesize{}0.579} & {\footnotesize{}0.731} &  & \tabularnewline\addlinespace[0.1cm]
\midrule 
\multicolumn{2}{c}{{\footnotesize{}DGP}} & {\footnotesize{}0} & {\footnotesize{}0} & {\footnotesize{}0} & {\footnotesize{}0} & {\footnotesize{}0} & {\footnotesize{}0} & {\footnotesize{}0} & {\footnotesize{}0} & {\footnotesize{}0} &  & {\footnotesize{}0}\tabularnewline\addlinespace[0.1cm]
\multicolumn{4}{l}{{\footnotesize{}Misspecified means}} &  &  &  &  &  &  &  &  & \tabularnewline
 & {\footnotesize{}{[}-2,-2{]}} & {\footnotesize{}975} & {\footnotesize{}3046} & {\footnotesize{}5854} & {\footnotesize{}8342} & {\footnotesize{}10464} & {\footnotesize{}12425} & {\footnotesize{}14161} & {\footnotesize{}15309} & {\footnotesize{}15071} &  & {\footnotesize{}21754}\tabularnewline
 & {\footnotesize{}{[}-0.4,-0.4{]}} & {\footnotesize{}141} & {\footnotesize{}269} & {\footnotesize{}374} & {\footnotesize{}476} & {\footnotesize{}532} & {\footnotesize{}565} & {\footnotesize{}597} & {\footnotesize{}548} & {\footnotesize{}433} &  & {\footnotesize{}817}\tabularnewline
 & {\footnotesize{}{[}0.4, 0.4{]}} & {\footnotesize{}204} & {\footnotesize{}347} & {\footnotesize{}431} & {\footnotesize{}501} & {\footnotesize{}555} & {\footnotesize{}547} & {\footnotesize{}526} & {\footnotesize{}437} & {\footnotesize{}288} &  & {\footnotesize{}721}\tabularnewline
 & {\footnotesize{}{[}2, 2{]}} & {\footnotesize{}8296} & {\footnotesize{}10002} & {\footnotesize{}10447} & {\footnotesize{}10145} & {\footnotesize{}9348} & {\footnotesize{}8151} & {\footnotesize{}6543} & {\footnotesize{}4638} & {\footnotesize{}2330} &  & {\footnotesize{}11033}\tabularnewline\addlinespace[0.1cm]
\multicolumn{5}{l}{{\footnotesize{}Misspecified variances}} &  &  &  &  &  &  &  & \tabularnewline
 & {\footnotesize{}{[}0.6,0.6{]}} & {\footnotesize{}144} & {\footnotesize{}137} & {\footnotesize{}102} & {\footnotesize{}67} & {\footnotesize{}39} & {\footnotesize{}24} & {\footnotesize{}27} & {\footnotesize{}51} & {\footnotesize{}103} &  & {\footnotesize{}156}\tabularnewline
 & {\footnotesize{}{[}0.8,0.8{]}} & {\footnotesize{}27} & {\footnotesize{}27} & {\footnotesize{}19} & {\footnotesize{}13} & {\footnotesize{}7} & {\footnotesize{}5} & {\footnotesize{}5} & {\footnotesize{}11} & {\footnotesize{}22} &  & {\footnotesize{}31}\tabularnewline
 & {\footnotesize{}{[}1.2,1.2{]}} & {\footnotesize{}18} & {\footnotesize{}18} & {\footnotesize{}15} & {\footnotesize{}9} & {\footnotesize{}5} & {\footnotesize{}3} & {\footnotesize{}4} & {\footnotesize{}9} & {\footnotesize{}16} &  & {\footnotesize{}22}\tabularnewline
 & {\footnotesize{}{[}1.4,1.4{]}} & {\footnotesize{}61} & {\footnotesize{}66} & {\footnotesize{}50} & {\footnotesize{}31} & {\footnotesize{}17} & {\footnotesize{}10} & {\footnotesize{}16} & {\footnotesize{}34} & {\footnotesize{}62} &  & {\footnotesize{}77}\tabularnewline\addlinespace[0.1cm]
\multicolumn{5}{l}{{\footnotesize{}Misspecified covariance}} &  &  &  &  &  &  &  & \tabularnewline
 & {\footnotesize{}-0.5} & {\footnotesize{}320} & {\footnotesize{}323} & {\footnotesize{}278} & {\footnotesize{}232} & {\footnotesize{}187} & {\footnotesize{}129} & {\footnotesize{}85} & {\footnotesize{}46} & {\footnotesize{}15} &  & {\footnotesize{}222}\tabularnewline
 & {\footnotesize{}0.1} & {\footnotesize{}31} & {\footnotesize{}39} & {\footnotesize{}40} & {\footnotesize{}34} & {\footnotesize{}31} & {\footnotesize{}27} & {\footnotesize{}19} & {\footnotesize{}13} & {\footnotesize{}5} &  & {\footnotesize{}37}\tabularnewline
 & {\footnotesize{}0.7} & {\footnotesize{}6} & {\footnotesize{}8} & {\footnotesize{}10} & {\footnotesize{}10} & {\footnotesize{}9} & {\footnotesize{}8} & {\footnotesize{}7} & {\footnotesize{}5} & {\footnotesize{}3} &  & {\footnotesize{}11}\tabularnewline
 & {\footnotesize{}0.9} & {\footnotesize{}19} & {\footnotesize{}30} & {\footnotesize{}35} & {\footnotesize{}40} & {\footnotesize{}38} & {\footnotesize{}38} & {\footnotesize{}35} & {\footnotesize{}27} & {\footnotesize{}16} &  & {\footnotesize{}49}\tabularnewline
\bottomrule
\end{tabular*}{\footnotesize\par}

{\small{}\vspace{-0.6cm}
}{\small\par}
\end{table}

\begin{table}[H]
{\footnotesize{}\caption{\linespread{0.55}\selectfont{}{\footnotesize{}For the simulated data,
comparison of candidate distributions using the LPM level set score
${\rm (LPMS^{'})^{\Gamma}}$ and $L^{2}$ score ${\rm LPMS^{'}}$
$(\times10^{5})$ for $k=1$. The scores were computed with $h$
chosen as the PDF of the uniform distribution on $[-3,3]^{2}$. Each
value is the score for a candidate distribution minus the score for
DGP. Lower values are better.}\label{tab:LPM level set_simulation}}
\vspace{0.1cm}
 }\linespread{0.45}\selectfont{}{\footnotesize{}}%
\begin{tabular*}{1\textwidth}{@{\extracolsep{\fill}}>{\raggedright}p{0.1cm}cccccccccccc}
\toprule 
 &  & \multicolumn{9}{c}{{\footnotesize{}${\rm (LPMS^{'})^{\Gamma}}$}} &  & {\footnotesize{}${\rm LPMS^{'}}$}\tabularnewline
\midrule 
 & $\alpha$ & {\footnotesize{}0.016} & {\footnotesize{}0.047} & {\footnotesize{}0.094} & {\footnotesize{}0.163} & {\footnotesize{}0.261} & {\footnotesize{}0.401} & {\footnotesize{}0.611} & {\footnotesize{}0.952} & {\footnotesize{}1.628} &  & \tabularnewline\addlinespace[0.1cm]
\midrule 
\multicolumn{2}{c}{{\footnotesize{}DGP}} & {\footnotesize{}0} & {\footnotesize{}0} & {\footnotesize{}0} & {\footnotesize{}0} & {\footnotesize{}0} & {\footnotesize{}0} & {\footnotesize{}0} & {\footnotesize{}0} & {\footnotesize{}0} &  & {\footnotesize{}0}\tabularnewline\addlinespace[0.1cm]
\multicolumn{4}{l}{{\footnotesize{}Misspecified means}} &  &  &  &  &  &  &  &  & \tabularnewline
 & {\footnotesize{}{[}-2,-2{]}} & {\footnotesize{}381} & {\footnotesize{}1588} & {\footnotesize{}3745} & {\footnotesize{}7050} & {\footnotesize{}11755} & {\footnotesize{}17914} & {\footnotesize{}25674} & {\footnotesize{}36623} & {\footnotesize{}55301} &  & {\footnotesize{}3278906}\tabularnewline
 & {\footnotesize{}{[}-0.4,-0.4{]}} & {\footnotesize{}99} & {\footnotesize{}232} & {\footnotesize{}408} & {\footnotesize{}624} & {\footnotesize{}865} & {\footnotesize{}1154} & {\footnotesize{}1538} & {\footnotesize{}2033} & {\footnotesize{}2669} &  & {\footnotesize{}52604}\tabularnewline
 & {\footnotesize{}{[}0.4, 0.4{]}} & {\footnotesize{}123} & {\footnotesize{}317} & {\footnotesize{}541} & {\footnotesize{}804} & {\footnotesize{}1118} & {\footnotesize{}1470} & {\footnotesize{}1864} & {\footnotesize{}2297} & {\footnotesize{}2820} &  & {\footnotesize{}32690}\tabularnewline
 & {\footnotesize{}{[}2, 2{]}} & {\footnotesize{}12288} & {\footnotesize{}20694} & {\footnotesize{}27612} & {\footnotesize{}33659} & {\footnotesize{}38675} & {\footnotesize{}42980} & {\footnotesize{}45236} & {\footnotesize{}45084} & {\footnotesize{}36098} &  & {\footnotesize{}261370}\tabularnewline\addlinespace[0.1cm]
\multicolumn{5}{l}{{\footnotesize{}Misspecified variances}} &  &  &  &  &  &  &  & \tabularnewline
 & {\footnotesize{}{[}0.6,0.6{]}} & {\footnotesize{}143} & {\footnotesize{}226} & {\footnotesize{}272} & {\footnotesize{}277} & {\footnotesize{}279} & {\footnotesize{}263} & {\footnotesize{}224} & {\footnotesize{}156} & {\footnotesize{}76} &  & {\footnotesize{}725}\tabularnewline
 & {\footnotesize{}{[}0.8,0.8{]}} & {\footnotesize{}19} & {\footnotesize{}36} & {\footnotesize{}50} & {\footnotesize{}57} & {\footnotesize{}57} & {\footnotesize{}55} & {\footnotesize{}52} & {\footnotesize{}36} & {\footnotesize{}21} &  & {\footnotesize{}178}\tabularnewline
 & {\footnotesize{}{[}1.2,1.2{]}} & {\footnotesize{}27} & {\footnotesize{}45} & {\footnotesize{}54} & {\footnotesize{}56} & {\footnotesize{}57} & {\footnotesize{}54} & {\footnotesize{}47} & {\footnotesize{}32} & {\footnotesize{}16} &  & {\footnotesize{}143}\tabularnewline
 & {\footnotesize{}{[}1.4,1.4{]}} & {\footnotesize{}76} & {\footnotesize{}132} & {\footnotesize{}170} & {\footnotesize{}189} & {\footnotesize{}199} & {\footnotesize{}193} & {\footnotesize{}172} & {\footnotesize{}133} & {\footnotesize{}72} &  & {\footnotesize{}577}\tabularnewline\addlinespace[0.1cm]
\multicolumn{5}{l}{{\footnotesize{}Misspecified covariance}} &  &  &  &  &  &  &  & \tabularnewline
 & {\footnotesize{}-0.5} & {\footnotesize{}560} & {\footnotesize{}877} & {\footnotesize{}1130} & {\footnotesize{}1351} & {\footnotesize{}1532} & {\footnotesize{}1714} & {\footnotesize{}1785} & {\footnotesize{}1821} & {\footnotesize{}1691} &  & {\footnotesize{}16643}\tabularnewline
 & {\footnotesize{}0.1} & {\footnotesize{}34} & {\footnotesize{}73} & {\footnotesize{}108} & {\footnotesize{}143} & {\footnotesize{}181} & {\footnotesize{}210} & {\footnotesize{}248} & {\footnotesize{}270} & {\footnotesize{}282} &  & {\footnotesize{}2835}\tabularnewline
 & {\footnotesize{}0.7} & {\footnotesize{}11} & {\footnotesize{}18} & {\footnotesize{}25} & {\footnotesize{}32} & {\footnotesize{}40} & {\footnotesize{}48} & {\footnotesize{}52} & {\footnotesize{}59} & {\footnotesize{}56} &  & {\footnotesize{}585}\tabularnewline
 & {\footnotesize{}0.9} & {\footnotesize{}27} & {\footnotesize{}50} & {\footnotesize{}80} & {\footnotesize{}112} & {\footnotesize{}140} & {\footnotesize{}173} & {\footnotesize{}202} & {\footnotesize{}241} & {\footnotesize{}241} &  & {\footnotesize{}2591}\tabularnewline
\bottomrule
\end{tabular*}{\footnotesize\par}

{\small{}\vspace{-0.6cm}
}{\small\par}
\end{table}

\begin{table}[H]
{\footnotesize{}\caption{\linespread{0.55}\selectfont{}{\footnotesize{}For the simulated data,
comparison of candidate distributions using the density level set
score ${\rm (DQS^{'})^{\Gamma}}$ and $L^{2}$ score ${\rm DQS^{'}}$
$(\times10^{5})$. The scores were computed with $h$ chosen as the
PDF of $\mathcal{N}(\mathbf{0},\mathbf{I}_{2})$. Each value is the
score for a candidate distribution minus the score for DGP. Lower
values are better.}\label{tab:CDF level set_simulation Gaussian}}
\vspace{0.1cm}
}\linespread{0.45}\selectfont{}{\footnotesize{} }%
\begin{tabular*}{1\textwidth}{@{\extracolsep{\fill}}>{\raggedright}p{0.1cm}cccccccccccc}
\toprule 
 &  & \multicolumn{9}{c}{{\footnotesize{}${\rm {\rm (DQS^{'})^{\Gamma}}}$}} &  & {\footnotesize{}${\rm DQS^{'}}$}\tabularnewline
\midrule 
 & {\footnotesize{}$\alpha$} & {\footnotesize{}0.018} & {\footnotesize{}0.037} & {\footnotesize{}0.055} & {\footnotesize{}0.073} & {\footnotesize{}0.092} & {\footnotesize{}0.110} & {\footnotesize{}0.129} & {\footnotesize{}0.147} & {\footnotesize{}0.165} &  & \tabularnewline\addlinespace[0.1cm]
\midrule 
\multicolumn{2}{c}{{\footnotesize{}DGP}} & {\footnotesize{}0} & {\footnotesize{}0} & {\footnotesize{}0} & {\footnotesize{}0} & {\footnotesize{}0} & {\footnotesize{}0} & {\footnotesize{}0} & {\footnotesize{}0} & {\footnotesize{}0} &  & {\footnotesize{}0}\tabularnewline\addlinespace[0.1cm]
\multicolumn{4}{l}{{\footnotesize{}Misspecified means}} &  &  &  &  &  &  &  &  & \tabularnewline
 & {\footnotesize{}{[}-2,-2{]}} & {\footnotesize{}4473} & {\footnotesize{}4483} & {\footnotesize{}3801} & {\footnotesize{}2976} & {\footnotesize{}2172} & {\footnotesize{}1454} & {\footnotesize{}862} & {\footnotesize{}431} & {\footnotesize{}119} &  & {\footnotesize{}803}\tabularnewline
 & {\footnotesize{}{[}-0.4,-0.4{]}} & {\footnotesize{}37} & {\footnotesize{}94} & {\footnotesize{}158} & {\footnotesize{}198} & {\footnotesize{}255} & {\footnotesize{}275} & {\footnotesize{}266} & {\footnotesize{}233} & {\footnotesize{}139} &  & {\footnotesize{}61}\tabularnewline
 & {\footnotesize{}{[}0.4, 0.4{]}} & {\footnotesize{}35} & {\footnotesize{}100} & {\footnotesize{}173} & {\footnotesize{}213} & {\footnotesize{}253} & {\footnotesize{}264} & {\footnotesize{}246} & {\footnotesize{}192} & {\footnotesize{}90} &  & {\footnotesize{}58}\tabularnewline
 & {\footnotesize{}{[}2, 2{]}} & {\footnotesize{}4472} & {\footnotesize{}4498} & {\footnotesize{}3809} & {\footnotesize{}2978} & {\footnotesize{}2171} & {\footnotesize{}1457} & {\footnotesize{}860} & {\footnotesize{}424} & {\footnotesize{}114} &  & {\footnotesize{}803}\tabularnewline\addlinespace[0.1cm]
\multicolumn{5}{l}{{\footnotesize{}Misspecified variances}} &  &  &  &  &  &  &  & \tabularnewline
 & {\footnotesize{}{[}0.6,0.6{]}} & {\footnotesize{}847} & {\footnotesize{}696} & {\footnotesize{}490} & {\footnotesize{}314} & {\footnotesize{}208} & {\footnotesize{}134} & {\footnotesize{}131} & {\footnotesize{}234} & {\footnotesize{}395} &  & {\footnotesize{}915}\tabularnewline
 & {\footnotesize{}{[}0.8,0.8{]}} & {\footnotesize{}67} & {\footnotesize{}54} & {\footnotesize{}37} & {\footnotesize{}6} & {\footnotesize{}10} & {\footnotesize{}21} & {\footnotesize{}65} & {\footnotesize{}134} & {\footnotesize{}172} &  & {\footnotesize{}70}\tabularnewline
 & {\footnotesize{}{[}1.2,1.2{]}} & {\footnotesize{}16} & {\footnotesize{}17} & {\footnotesize{}10} & {\footnotesize{}3} & {\footnotesize{}25} & {\footnotesize{}96} & {\footnotesize{}211} & {\footnotesize{}357} & {\footnotesize{}85} &  & {\footnotesize{}29}\tabularnewline
 & {\footnotesize{}{[}1.4,1.4{]}} & {\footnotesize{}42} & {\footnotesize{}39} & {\footnotesize{}23} & {\footnotesize{}38} & {\footnotesize{}195} & {\footnotesize{}583} & {\footnotesize{}751} & {\footnotesize{}357} & {\footnotesize{}85} &  & {\footnotesize{}80}\tabularnewline\addlinespace[0.1cm]
\multicolumn{5}{l}{{\footnotesize{}Misspecified covariance}} &  &  &  &  &  &  &  & \tabularnewline
 & {\footnotesize{}-0.5} & {\footnotesize{}342} & {\footnotesize{}540} & {\footnotesize{}644} & {\footnotesize{}641} & {\footnotesize{}596} & {\footnotesize{}474} & {\footnotesize{}324} & {\footnotesize{}191} & {\footnotesize{}38} &  & {\footnotesize{}141}\tabularnewline
 & {\footnotesize{}0.1} & {\footnotesize{}61} & {\footnotesize{}105} & {\footnotesize{}138} & {\footnotesize{}119} & {\footnotesize{}116} & {\footnotesize{}77} & {\footnotesize{}94} & {\footnotesize{}129} & {\footnotesize{}85} &  & {\footnotesize{}35}\tabularnewline
 & {\footnotesize{}0.7} & {\footnotesize{}79} & {\footnotesize{}79} & {\footnotesize{}68} & {\footnotesize{}36} & {\footnotesize{}43} & {\footnotesize{}38} & {\footnotesize{}64} & {\footnotesize{}62} & {\footnotesize{}75} &  & {\footnotesize{}33}\tabularnewline
 & {\footnotesize{}0.9} & {\footnotesize{}920} & {\footnotesize{}791} & {\footnotesize{}607} & {\footnotesize{}436} & {\footnotesize{}320} & {\footnotesize{}234} & {\footnotesize{}217} & {\footnotesize{}288} & {\footnotesize{}389} &  & {\footnotesize{}494}\tabularnewline
\bottomrule
\end{tabular*}{\footnotesize\par}

{\small{}\vspace{-0.6cm}
}{\small\par}
\end{table}

\begin{table}[H]
{\footnotesize{}\caption{\linespread{0.55}\selectfont{}{\footnotesize{}For the simulated data,
comparison of candidate distributions using the density level set
score ${\rm (DQS^{'})^{\Gamma}}$ and $L^{2}$ score ${\rm DQS^{'}}$
$(\times10^{5})$. The scores were computed with $h$ chosen as the
PDF of $\mathcal{N}(\mathbf{0},2\times\mathbf{I}_{2})$. Each value
is the score for a candidate distribution minus the score for DGP.
Lower values are better.}\label{tab:CDF level set_simulation Gaussian 2}}
\vspace{0.1cm}
}\linespread{0.45}\selectfont{}{\footnotesize{} }%
\begin{tabular*}{1\textwidth}{@{\extracolsep{\fill}}>{\raggedright}p{0.1cm}cccccccccccc}
\toprule 
 &  & \multicolumn{9}{c}{{\footnotesize{}${\rm {\rm (DQS^{'})^{\Gamma}}}$}} &  & {\footnotesize{}${\rm DQS^{'}}$}\tabularnewline
\midrule 
 & {\footnotesize{}$\alpha$} & {\footnotesize{}0.018} & {\footnotesize{}0.037} & {\footnotesize{}0.055} & {\footnotesize{}0.073} & {\footnotesize{}0.092} & {\footnotesize{}0.110} & {\footnotesize{}0.129} & {\footnotesize{}0.147} & {\footnotesize{}0.165} &  & \tabularnewline\addlinespace[0.1cm]
\midrule 
\multicolumn{2}{c}{{\footnotesize{}DGP}} & {\footnotesize{}0} & {\footnotesize{}0} & {\footnotesize{}0} & {\footnotesize{}0} & {\footnotesize{}0} & {\footnotesize{}0} & {\footnotesize{}0} & {\footnotesize{}0} & {\footnotesize{}0} &  & {\footnotesize{}0}\tabularnewline\addlinespace[0.1cm]
\multicolumn{4}{l}{{\footnotesize{}Misspecified means}} &  &  &  &  &  &  &  &  & \tabularnewline
 & {\footnotesize{}{[}-2,-2{]}} & {\footnotesize{}2751} & {\footnotesize{}2688} & {\footnotesize{}2300} & {\footnotesize{}1826} & {\footnotesize{}1360} & {\footnotesize{}960} & {\footnotesize{}621} & {\footnotesize{}342} & {\footnotesize{}116} &  & {\footnotesize{}506}\tabularnewline
 & {\footnotesize{}{[}-0.4,-0.4{]}} & {\footnotesize{}58} & {\footnotesize{}103} & {\footnotesize{}150} & {\footnotesize{}162} & {\footnotesize{}169} & {\footnotesize{}165} & {\footnotesize{}163} & {\footnotesize{}118} & {\footnotesize{}33} &  & {\footnotesize{}42}\tabularnewline
 & {\footnotesize{}{[}0.4, 0.4{]}} & {\footnotesize{}57} & {\footnotesize{}115} & {\footnotesize{}166} & {\footnotesize{}173} & {\footnotesize{}170} & {\footnotesize{}168} & {\footnotesize{}159} & {\footnotesize{}122} & {\footnotesize{}60} &  & {\footnotesize{}44}\tabularnewline
 & {\footnotesize{}{[}2, 2{]}} & {\footnotesize{}2756} & {\footnotesize{}2703} & {\footnotesize{}2316} & {\footnotesize{}1844} & {\footnotesize{}1381} & {\footnotesize{}977} & {\footnotesize{}636} & {\footnotesize{}351} & {\footnotesize{}121} &  & {\footnotesize{}511}\tabularnewline\addlinespace[0.1cm]
\multicolumn{5}{l}{{\footnotesize{}Misspecified variances}} &  &  &  &  &  &  &  & \tabularnewline
 & {\footnotesize{}{[}0.6,0.6{]}} & {\footnotesize{}639} & {\footnotesize{}484} & {\footnotesize{}333} & {\footnotesize{}200} & {\footnotesize{}113} & {\footnotesize{}86} & {\footnotesize{}106} & {\footnotesize{}164} & {\footnotesize{}251} &  & {\footnotesize{}527}\tabularnewline
 & {\footnotesize{}{[}0.8,0.8{]}} & {\footnotesize{}60} & {\footnotesize{}44} & {\footnotesize{}38} & {\footnotesize{}3} & {\footnotesize{}11} & {\footnotesize{}18} & {\footnotesize{}50} & {\footnotesize{}80} & {\footnotesize{}107} &  & {\footnotesize{}44}\tabularnewline
 & {\footnotesize{}{[}1.2,1.2{]}} & {\footnotesize{}20} & {\footnotesize{}14} & {\footnotesize{}9} & {\footnotesize{}8} & {\footnotesize{}21} & {\footnotesize{}40} & {\footnotesize{}112} & {\footnotesize{}164} & {\footnotesize{}24} &  & {\footnotesize{}15}\tabularnewline
 & {\footnotesize{}{[}1.4,1.4{]}} & {\footnotesize{}53} & {\footnotesize{}34} & {\footnotesize{}23} & {\footnotesize{}35} & {\footnotesize{}118} & {\footnotesize{}314} & {\footnotesize{}376} & {\footnotesize{}164} & {\footnotesize{}24} &  & {\footnotesize{}43}\tabularnewline\addlinespace[0.1cm]
\multicolumn{5}{l}{{\footnotesize{}Misspecified covariance}} &  &  &  &  &  &  &  & \tabularnewline
 & {\footnotesize{}-0.5} & {\footnotesize{}419} & {\footnotesize{}526} & {\footnotesize{}546} & {\footnotesize{}476} & {\footnotesize{}382} & {\footnotesize{}285} & {\footnotesize{}198} & {\footnotesize{}96} & {\footnotesize{}16} &  & {\footnotesize{}112}\tabularnewline
 & {\footnotesize{}0.1} & {\footnotesize{}82} & {\footnotesize{}101} & {\footnotesize{}113} & {\footnotesize{}94} & {\footnotesize{}69} & {\footnotesize{}48} & {\footnotesize{}48} & {\footnotesize{}61} & {\footnotesize{}24} &  & {\footnotesize{}24}\tabularnewline
 & {\footnotesize{}0.7} & {\footnotesize{}69} & {\footnotesize{}64} & {\footnotesize{}62} & {\footnotesize{}30} & {\footnotesize{}27} & {\footnotesize{}16} & {\footnotesize{}41} & {\footnotesize{}38} & {\footnotesize{}42} &  & {\footnotesize{}24}\tabularnewline
 & {\footnotesize{}0.9} & {\footnotesize{}681} & {\footnotesize{}554} & {\footnotesize{}418} & {\footnotesize{}295} & {\footnotesize{}214} & {\footnotesize{}171} & {\footnotesize{}173} & {\footnotesize{}217} & {\footnotesize{}280} &  & {\footnotesize{}307}\tabularnewline
\bottomrule
\end{tabular*}{\footnotesize\par}

{\small{}\vspace{-0.6cm}
}{\small\par}
\end{table}

\medskip{}

{\color{\mycolor}

\section{Practical Applications}

In this section, we apply the $L^{2}$ scoring rules
and level set scores to two practical applications. In Section 6.1,
we use the $L^{2}$ scoring rules to estimate weights for combining
distributional forecasts, and in Section 6.2, we use the CDF level
set score to estimate conditional Value-at-Risk (CoVaR). 

\subsection{Using $L^{2}$ Scores to Combine Distributional
Forecasts\label{subsec:Combining L2}}

Combining is a popular and pragmatic way to improve
forecast accuracy. In recent years, interest has increased in methods
for combining probabilistic forecasts \citep{aastveit2018evolution,winkler2019probability}.
We show that, for a combination of distributional forecasts, the $L^{2}$
scoring rules enable us to estimate the combining weights as a quadratic
optimization algorithm, which is both theoretically and computationally
appealing. Let $\left\{ P_{\mathbf{X}_{it}}\right\} _{i=1}^{N}$ be
$N$ distributional forecasts for $\mathbf{Y}_{t}$, where $\mathbf{X}_{it}$ has density
$ f_{\mathbf{X}_{it}}$ for $t = 1, 2, \dots, T$. The most common
method for combining distributional forecasts is to produce a mixture
distribution of the form $\sum_{i=1}^{N}\eta_{i}P_{\mathbf{X}_{it}}$,
which has PDF $\sum_{i=1}^{N}\eta_{i}f_{\mathbf{X}_{it}}$, where
$\sum_{i=1}^{N}\eta_{i}=1$ and $\eta_{i}\geq0$. The standard approach
to estimating the combining weights $\eta_{i}$ is to minimize the
log score, which can be viewed as maximum likelihood estimation (see, for example,
\citealt{hall2007combining,gneiting2013combining}). Our alternative
proposal is to estimate the weights by minimizing an $L^{2}$ scoring
rule.

Let us consider the form of the $L^{2}$ scores in
(\ref{eq: L2 score equivalent renormalization physical}). If we expand the square bracket and rearrange the terms in
(\ref{eq: L2 score equivalent renormalization physical}), we find that, with $T$ in-sample observations, the average $L^{2}$ score, $\frac{1}{T}\sum_{t=1}^{T}S^{'}\left(\sum_{i=1}^{N}\eta_{i}P_{\mathbf{X}_{it}},\mathbf{y}_{t};w, h \right)$,
can be simplified into a quadratic function of the weights
$\eta_{i}$. In fact, it can be equivalently written, as $\bm{\eta}^{'}B\mathbf{\bm{\eta}}+\mathbf{c^{'}\mathbf{\mathbf{\bm{\eta}}}}$, where $\mathbf{\bm{\eta}} = [\eta_1, \eta_2, \dots, \eta_N]^{'}$, $B$ is an $N \times N$ positive semi-definite matrix with $B_{i,j}=\frac{1}{T}\sum_{t=1}^{T}\intop_{\mathbb{R}^{d}}(f_{\mathbf{X}_{it}}*w)(\mathbf{z})(f_{\mathbf{X}_{jt}}*w)(\mathbf{z})h(\mathbf{z})d\mathbf{z}$ and $\mathbf{c} = [c_1, c_2, \dots, c_N]^{'}$ with $c_{i}=-\frac{2}{T}\sum_{t=1}^{T}\intop_{\mathbb{R}^{d}}\left(f_{\mathbf{X}_{it}}*w\right)(\mathbf{z})w(\mathbf{z}-\mathbf{y})h(\mathbf{z})d\mathbf{z}$. Hence, to find the optimal weights, we simply need to solve a quadratic optimization problem with linear constraints $\mathbf{1}_{N}^{'}\bm{\eta}=1$
and $\bm{\eta}\geq0$, where $\mathbf{1}_{N}$ is a column vector
of ones.

Our proposed approach has several desirable properties
when compared with estimating the combining weights by minimizing
the log score. First, the $L^{2}$ scores offer greater flexibility.
This is because, unlike the log score, the $L^{2}$ scores allow the
PDFs to have zero values. Furthermore, by contrast with the log score,
the CDF-based $\text{MCRPS}$ and $\text{MCRPS}{}^{'}$, can be used
for estimation in applications where PDFs are not available, such
as when the distributional forecast is produced via simulation (see, for example, \citealt{taylor2017probabilistic}).
Second, the quadratic optimization problem is computationally efficient,
as it can be solved numerically in polynomial time \citep{ye1989extension}.
Third, the matrix $B$ reveals useful information regarding each individual
forecast. Any individual forecast that is perfectly correlated with
the others in the sense of the matrix $B$ can be viewed as redundant,
as it has no impact on the minimal value of the $L^{2}$ score. By
removing the corresponding rows (columns), we can assume that the
matrix $B$ is positive definite, which will induce a unique solution for the optimal combining weights.  Finally, our proposed approach generalizes the seminal
work of \citet{makridakis1983averages} on point forecasting. Hence,
recent developments for combining point forecasts (see, for example,
\citealt{soule2021heuristic}) can be incorporated in our approach
for combining probabilistic forecasts. Our work also includes, as a special case, the work
of \citet{hora2015calibration}, who use the quadratic score to estimate
combining weights.

To illustrate our proposal, we use forecasts provided
by multiple experts in the ECB Survey of Professional Forecasters.
The data consists of rolling year-ahead predictions for quarterly
inflation, GDP and unemployment in the euro area for the 91 quarters
up to the third quarter of 2021. The distributional forecasts are
submitted in the form of probabilities for the economic variable falling
in pre-specified bins. There were many cases where forecasters submitted
zero probability for at least one bin, implying that the PDFs took
zero values. This meant that the log score could not be directly used as the
basis for estimating combining weights. Although there were more than
100 individual forecasters, many did not provide forecasts for every
quarter. We kept only the individual forecasters for which forecasts
were missing for eight or fewer quarters. This left six forecasters for GDP, seven for inflation, and five for unemployment. For these forecasters,
any missing distributional forecasts were replaced by the average
of all the other available forecasts for that quarter.

We considered three combining methods: the simple
average, and the proposed combining approach with weights estimated
using the ${\rm DQS^{'}}$ and ${\rm CRPS^{'}}$. To compute the $L^{2}$
scores, for each series, $M = 10^{4}$ values were sampled from $h$,
which we chose to be the uniform distribution between the minimum
and maximum values for the entire period. We repeatedly re-estimated
model parameters using a rolling window of 48 quarters, which delivered
43 out-of-sample forecasts. We also used ${\rm DQS^{'}}$ and ${\rm CRPS^{'}}$
for evaluation. We present the results in Table 6, where the scores
corresponding to the benchmark simple average,
which is used as the reference, have been subtracted
from each candidate. It can be seen that all the scores associated
with ${\rm DQS^{'}}$ and ${\rm CRPS^{'}}$ are negative, indicating
that they outperformed the simple average.

\begin{table}[H]
{\color{\mycolor} {}
\footnotesize{}\caption{\linespread{0.5}\selectfont{}{\footnotesize{} Comparison
of combining methods. Each value is a score  for a method minus the
score for the simple average $(\times10^{3})$. Lower values are better. $^{\dag}$
and {*} indicate significantly less than zero at 10\% and 5\% significance
levels}.} \label{tab:Density level set_TS}
\smallskip{}
}{\footnotesize\par}

\linespread{0.55}\selectfont{}\textcolor{\mycolor}{\footnotesize{}}%
\begin{tabular*}{1\textwidth}{@{\extracolsep{\fill}}lllccccccccccc}
\toprule 
 &  &  & \multicolumn{3}{c}{\textcolor{\mycolor}{\footnotesize{}GDP}} &  & \multicolumn{3}{c}{\textcolor{\mycolor}{\footnotesize{}Inflation}} &  & \multicolumn{3}{c}{\textcolor{\mycolor}{\footnotesize{}Unemployment}}\tabularnewline
\midrule 
 &  &  & \multicolumn{3}{c}{\textcolor{\mycolor}{\footnotesize{}Combining method}} &  & \multicolumn{3}{c}{\textcolor{\mycolor}{\footnotesize{}Combining method}} &  & \multicolumn{3}{c}{\textcolor{\mycolor}{\footnotesize{}Combining method}}\tabularnewline
\multicolumn{2}{l}{} &  & \textcolor{\mycolor}{\footnotesize{}Simple Avg} & \textcolor{\mycolor}{\footnotesize{}${\rm DQS^{'}}$} & \textcolor{\mycolor}{\footnotesize{}${\rm CRPS^{'}}$} &  & \textcolor{\mycolor}{\footnotesize{}Simple Avg} & \textcolor{\mycolor}{\footnotesize{}${\rm DQS^{'}}$} & \textcolor{\mycolor}{\footnotesize{}${\rm CRPS^{'}}$} &  & \textcolor{\mycolor}{\footnotesize{}Simple Avg} & \textcolor{\mycolor}{\footnotesize{}${\rm DQS^{'}}$} & \textcolor{\mycolor}{\footnotesize{}${\rm CRPS^{'}}$}\tabularnewline
\midrule
 & \textcolor{\mycolor}{\footnotesize{}${\rm DQS^{'}}$} &  & \textcolor{\mycolor}{\footnotesize{}0} & \textcolor{\mycolor}{\footnotesize{}-42} & \textcolor{\mycolor}{\footnotesize{}$-97^{*}$} &  & \textcolor{\mycolor}{\footnotesize{}0} & \textcolor{\mycolor}{\footnotesize{}$-128$} & \textcolor{\mycolor}{\footnotesize{}$-258^{*}$} &  & \textcolor{\mycolor}{\footnotesize{}0} & \textcolor{\mycolor}{\footnotesize{}$-630^{*}$} & \textcolor{\mycolor}{\footnotesize{}$-651^{*}$}\tabularnewline
 & \textcolor{\mycolor}{\footnotesize{}${\rm CRPS^{'}}$} &  & \textcolor{\mycolor}{\footnotesize{}0} & \textcolor{\mycolor}{\footnotesize{}$-42^{*}$} & \textcolor{\mycolor}{\footnotesize{}$-5$} &  & \textcolor{\mycolor}{\footnotesize{}0} & \textcolor{\mycolor}{\footnotesize{}$-164^{*}$} & \textcolor{\mycolor}{\footnotesize{}$-31^{\dag}$} &  & \textcolor{\mycolor}{\footnotesize{}0} & \textcolor{\mycolor}{\footnotesize{}$-290^{*}$} & \textcolor{\mycolor}{\footnotesize{}$-301^{*}$}\tabularnewline
\bottomrule
\end{tabular*}{\footnotesize\par}

\textcolor{\mycolor}{\small{}\vspace{-0.6cm}
}{\small\par}
\end{table}

\subsection{Using the CDF Level Set Score for CoVaR Estimation}

CoVaR is a popular measure of systemic risk, which
assesses the probability that a market incurs a heavy loss given that
an individual stock is already in distress (see, for example, \citealt{tobias2016covar}, \citealt{Dimitriadis2022}). Let $Y_{1t}$ and $Y_{2t}$ denote returns for a market index
and an individual stock, respectively. The CoVaR of $Y_{1t}$ conditional
on $Y_{2t}$, at probability level $\theta$, is defined as the scalar,
$\text{CoVaR}_{Y_{1t}|Y_{2t}}$, that satisfies $\mathbb{P}\left(Y_{1t}<\text{CoVaR}_{Y_{1t}|Y_{2t}}|Y_{2t}<q_{Y_{2t}}(\theta)\right)=\theta$,
where $q_{Y_{2t}}(\theta)$ denotes the $\theta$ quantile of $Y_{2t}$.
In this paper, to forecast $\text{CoVaR}_{Y_{1t}|Y_{2t}}$, we consider
a standard approach in the literature, which is to use a GARCH-Copula
model (see, for example, \citealt{patton2012review}). The approach
involves two steps. First, marginal distributions of $Y_{1t}$ and
$Y_{2t}$ are estimated using GARCH models, and then a copula is used
to capture the dependence between the marginal distributions. 

We note that ${\color{\mycolor}\mathbb{P}\left(Y_{1t}<\text{CoVaR}_{Y_{1t}|Y_{2t}},Y_{2t}<q_{Y_{2t}}(\theta)\right)=\theta^{2}}$.
This indicates that the pair $\left(\text{CoVaR}_{Y_{1t}|Y_{2t}},q_{Y_{2t}}(\theta)\right)$
is precisely on the $\theta^{2}$ CDF level set of $(Y_{1t},Y_{2t})$. This prompts
us to propose the estimation of the copula parameters by minimizing
the score for the $\theta^{2}$ CDF level set of $(Y_{1t},Y_{2t})$. We argue that
this will deliver better forecasting accuracy, in comparison with
the standard maximum likelihood approach, because the CDF level set
score enables improved fit in the tails.

To illustrate the proposed CoVaR estimation method
based on the CDF level set score, we set $Y_{1t}$ as the S\&P 500
index, and for $Y_{2t}$, we consider, in turn, the three stocks of the S\&P
500 that were analyzed by \citet{diks2020comparing}: Alcoa (AA),
MacDonald's (MCD), and Merck (MRK). We used the 5000 daily log returns
recorded between 21 February 2002 and 31 December 2021. Using a rolling
window of 2000 observations, we repeatedly re-estimated model parameters
to generate 3000 out-of-sample day-ahead forecasts for the joint distribution.
We used a GARCH(1,1) model with Student-t distribution for the marginal
distributions, and a Gaussian copula to model the dependence.
The CoVaR probability level $\theta$ was set to $10\%$. We estimated
the copula parameters first using maximum likelihood estimation, and then using
our method based on the score for the $\theta^{2}$=1\% CDF level set. This CDF level set score was computed using $10^{4}$ values sampled from the PDF $h$, which we set as the uniform distribution on {$[-0.3,0.3]^2$}.

Table~\ref{tab:CoVaR results} compares the CoVaR
forecasting performance of the two approaches using three measures.
First, we computed the CoVaR coverage as the proportion of the days
for which both of the following were true: $Y_{2t}<q_{Y_{2t}}(\theta)$
and $Y_{1t}<\text{CoVaR}_{Y_{1t}|Y_{2t}}$. In the table, the coverage
values for the proposed method are closer to the desired value of
1\% (=$\theta^{2}$) compared to the maximum likelihood approach. Second,
as $\text{CoVaR}_{Y_{1t}|Y_{2t}}$ is the $\theta$ quantile of $Y_{1t}$
given $Y_{2t}<q_{Y_{2t}}(\theta)$, we report the $\theta$
quantile score for $\text{CoVaR}_{Y_{1t}|Y_{2t}}$ for periods in
which $Y_{2t}<q_{Y_{2t}}(\theta)$. Third, we report
the 1\% CDF level set score.
For both these scores, lower values are better. For clarity, the scores
corresponding to the benchmark maximum likelihood approach are
chosen as the reference and subtracted from those
for the proposed method. In all three rows,
the proposed method has negative values for the scores, indicating it outperforms the benchmark in estimating the CoVaR for all three stocks.
\begin{table}[H]
\textcolor{\mycolor}{\footnotesize{}\caption{\linespread{0.5}\selectfont{}\textcolor{\mycolor}{\footnotesize{}Comparison
of CoVaR forecasting methods. Coverage closer to 1\% is preferred.
For the CoVaR quantile and level set scores, each value is the score
for the proposed approach minus the score for the maximum likelihood
benchmark $(\times10^{7})$. Lower values are better. $^{\dag}$
and {*} indicate significantly less than zero at 10\% and 5\% significance
levels.}\textcolor{purple}{\label{tab:CoVaR results}}}
\smallskip{}
}{\footnotesize\par}

\linespread{0.55}\selectfont{}\textcolor{\mycolor}{\footnotesize{}}%
\begin{tabular*}{1\textwidth}{@{\extracolsep{\fill}}lllcccccccc}
\hline 
 &  &  & \multicolumn{2}{c}{\textcolor{\mycolor}{\footnotesize{}Coverage (\%)}} &  & \multicolumn{2}{c}{\textcolor{\mycolor}{\footnotesize{}CoVaR quantile score}} &  & \multicolumn{2}{c}{\textcolor{\mycolor}{\footnotesize{}Level set score}}\tabularnewline
\hline 
 &  &  & \textcolor{\mycolor}{\footnotesize{}Max likelihood} & \textcolor{\mycolor}{\footnotesize{}Proposed} &  & \textcolor{\mycolor}{\footnotesize{}Max likelihood} & \textcolor{\mycolor}{\footnotesize{}Proposed} &  & \textcolor{\mycolor}{\footnotesize{}Max likelihood} & \textcolor{\mycolor}{\footnotesize{}Proposed}\tabularnewline
\multicolumn{2}{l}{\textcolor{\mycolor}{\footnotesize{}AA \& SP500}} &  & \textcolor{\mycolor}{\footnotesize{}1.67} & \textcolor{\mycolor}{\footnotesize{}1.37} &  & \textcolor{\mycolor}{\footnotesize{}0.00} & \textcolor{\mycolor}{\footnotesize{}$-6.04$} &  & \textcolor{\mycolor}{\footnotesize{}0.00} & \textcolor{\mycolor}{\footnotesize{}$-13.92$}\tabularnewline
\multicolumn{2}{l}{\textcolor{\mycolor}{\footnotesize{}MRK \& SP500}} &  & \textcolor{\mycolor}{\footnotesize{}1.77} & \textcolor{\mycolor}{\footnotesize{}1.37} &  & \textcolor{\mycolor}{\footnotesize{}0.00} & \textcolor{\mycolor}{\footnotesize{}$-9.88^{\dag}$} &  & \textcolor{\mycolor}{\footnotesize{}0.00} & \textcolor{\mycolor}{\footnotesize{}$-40.29$}\tabularnewline
\multicolumn{2}{l}{\textcolor{\mycolor}{\footnotesize{}MCD \& SP500}} &  & \textcolor{\mycolor}{\footnotesize{}1.60} & \textcolor{\mycolor}{\footnotesize{}1.27} &  & \textcolor{\mycolor}{\footnotesize{}0.00} & \textcolor{\mycolor}{\footnotesize{}$-8.10^{*}$} &  & \textcolor{\mycolor}{\footnotesize{}0.00} & \textcolor{\mycolor}{\footnotesize{}$-12.10$}\tabularnewline
\hline 
\end{tabular*}{\footnotesize\par}

\textcolor{\mycolor}{\small{}\vspace{-0.6cm}
}{\small\par}
\end{table}}

\section{Conclusion\label{sec:Conclusion}}

Forecasts of multivariate distributions and level sets are needed
to support decision making in a variety of contexts. In this paper,
we have studied the scoring rules for multivariate distributions and
scoring functions for level sets. The paper has several novel contributions.
Firstly, we propose the class of $L^{2}$ scoring rules for multivariate
distributions, for which the existing quadratic score and MCRPS are
specific examples. The $L^{2}$ scoring functions can easily generate
new scoring rules for multivariate distributions, and we demonstrate
this with the introduction of the LPMS, a new scoring rule based on
the lower partial moments. Secondly, by decomposing the $L^{2}$ scoring
rules, we obtain a unified approach for generating scoring functions
for level sets, including the scoring functions for density level
sets, CDF level sets, and LPM level sets. Thirdly, we propose a simple
numerical approach for computing the $L^{2}$ scoring rules and the
scoring functions for level sets. \textcolor{\mycolor}{Finally, we performed
a simulation study to provide support for the theoretical properties
of our new scores, and we used real data to illustrate their practical
usefulness for forecast combining and CoVaR estimation}.

\ACKNOWLEDGMENT{The authors are grateful to the Area Editor, the Associate Editor and three reviewers for providing very helpful comments. The authors are also grateful to Tobias Fissler and participants at the INFORMS Advances in Decision Analysis conference held at Bocconi University in Milan in 2019 for insightful feedback. }

\begin{SingleSpacedXII}{\small{}\bibliographystyle{chicago}
\bibliography{Reference_all}
}\end{SingleSpacedXII}


\begin{APPENDIX}{Assumption in Theorem \ref{thm:main}}
\begin{assumption}\label{assumption: one-1} 

\begin{enumerate}[label=(\alph*)]
    \item \textcolor{\mycolor}{Let $\lambda$ be a Borel measure on $\R^{d}$
whose Radon--Nikodym derivative with respect to the Lebesgue measure
is a non-negative integrable function $h$, i.e., $d\lambda(\mathbf{z})=h(\mathbf{z})d\mathbf{z}$.
Let $w$ be a local Borel measure such that $f_{{\bf Y}}\ast w$ is
a well-defined local Borel measure $\forall P_{{\bf Y}}\in\mathcal{V}^{d}$.}

    \item \textcolor{\mycolor}{If $f_{{\bf Y}}\ast w$ and $w$ are square-integrable
with respect to $\lambda$, and $\int_{\R^{d}}w^{2}(\mathbf{z}-\bullet)h({\bf z})f_{\mathbf{Y}}(\bullet)\,d{\bf z}$
is integrable with respect to the PDF $f_{\mathbf{Y}}$, then (\ref{eq: L2 Score physical})
is well-defined and finite.}

    \item \textcolor{\mycolor}{ If $\left(f_{\mathbf{X}}*w\right)(\bullet)w(\mathbf{\bullet}-\mathbf{y})$
is integrable with respect to $\lambda$, and $\mathbf{\int}_{\R^{d}}\left(f_{\mathbf{X}}*w\right)(\mathbf{z})w(\mathbf{z}-\mathbf{\bullet})h({\bf z})f_{\mathbf{Y}}(\mathbf{\bullet})\,d{\bf z}$
is integrable for any distribution $P_{{\bf Y}}\in\mathcal{V}^{d}$,
then (\ref{eq: L2 score equivalent renormalization physical}) is
well-defined and finite.}
    
\end{enumerate}

\end{assumption}

\end{APPENDIX}

\begin{APPENDIX}{Proofs of Theorems}

\proof{Proof of Theorem \ref{thm:main}} Consider the divergence $\Delta:=\mathbb{E}_{P_{\mathbf{Y}}}\left[S(P_{{\bf X}},\mathbf{\bullet};w)\right]-\mathbb{E}_{P_{\mathbf{Y}}}\left[S(P_{{\bf Y}},\mathbf{\bullet};w)\right].$
One may compute it as follows: 
\begin{align*}
\Delta & =\intop_{\mathbb{R}^{d}}\intop_{\mathbb{R}^{d}}\Big[\left(f_{\mathbf{X}}*w\right)(\mathbf{z})-w(\mathbf{z}-\mathbf{s})\Big]^{2}f_{\mathbf{Y}}(\mathbf{s})h({\bf z})d{\bf z}d\mathbf{s}-\intop_{\mathbb{R}^{d}}\intop_{\mathbb{R}^{d}}\Big[\left(f_{\mathbf{Y}}*w\right)(\mathbf{z})-w(\mathbf{z}-\mathbf{s})\Big]^{2}f_{\mathbf{Y}}(\mathbf{s})h({\bf z})d{\bf z}d\mathbf{s}\\
 & =\intop_{\mathbb{R}^{d}}\intop_{\mathbb{R}^{d}}\left\{ \left(f_{\mathbf{X}}*w\right)^{2}(\mathbf{z})-\left(f_{\mathbf{Y}}*w\right)^{2}(\mathbf{z})-2\Big[\left(f_{\mathbf{X}}*w\right)(\mathbf{z})-\left(f_{\mathbf{Y}}*w\right)(\mathbf{z})\Big]w(\mathbf{z}-\mathbf{s})\right\} f_{\mathbf{Y}}(\mathbf{s})\,d\mathbf{s}\,h({\bf z})\,d{\bf z} \\
 & =\intop_{\mathbb{R}^{d}}\left\{ \left(f_{\mathbf{X}}*w\right)^{2}(\mathbf{z})-\left(f_{\mathbf{Y}}*w\right)^{2}(\mathbf{z})-2\Big[\left(f_{\mathbf{X}}*w\right)(\mathbf{z})-\left(f_{\mathbf{Y}}*w\right)(\mathbf{z})\Big]\left(f_{\mathbf{Y}}*w\right)(\mathbf{z})\right\} h({\bf z})\,d{\bf z}\\
 & =\intop_{\mathbb{R}^{d}}\Big\{ f_{\mathbf{X}}*w-f_{\mathbf{Y}}*w\Big\}^{2}(\mathbf{z})h({\bf z})\,d{\bf z}\geq0,
\end{align*}
where the penultimate equality holds by Fubini's theorem and the definition
of convolution.

For strict properness, we assume $\Delta=0$. If $h$ is non-zero
$\mathcal{L}^{d}$-a.e., we obtain that $f_{\mathbf{X}}*w=f_{\mathbf{Y}}*w$
$\mathcal{L}^{d}$-a.e.. If $f_{\mathbf{Y}}*w$ uniquely characterizes
$P_{{\bf Y}}$ $\forall P_{{\bf Y}}\in\mathcal{V}^{d}$, we obtain
$P_{\mathbf{X}}=P_{\mathbf{Y}}$. \Halmos

\endproof
\proof{Derivation and proof of Theorem \ref{thm:scoring function for level set}}

 {\color{\mycolor}
(i) Derivation of (\ref{eq:scoring functions for level sets})
in Theorem \ref{thm:scoring function for level set}. The layer
cake representation states that if $\mu$ is a Borel measure on $\R^{d}$
and $g:\R^{d}\map[0,\infty)$ is a $\mu$-measurable function, then for
any $p\in[1,\infty)$, there holds
\[
\intop_{\mathbb{R}^{d}}g(\mathbf{z})^{p}\,d\mu(\mathbf{z})=\intop_{0}^{\infty}p\alpha^{p-1}{\mu}\big\{ g\geq\alpha\big\}\,d\alpha=\intop_{0}^{\infty}p\alpha^{p-1}{\mu}\left\{ L(g;\alpha)\right\} \,d\alpha,
\]
where $\mu\{\mathbf{\bullet}\}$ denotes the measure of a set under
$\mu$. 

Let us define two Borel measures, $\lambda$ as $d\lambda(\mathbf{z}):=h(\mathbf{z})\,d\mathbf{z}$
, and $\mu_{\yy,w}$ as $d\mu_{\yy,w}(\mathbf{z}):=w(\mathbf{z}-\mathbf{y})h(\mathbf{z})d\mathbf{z}$. Then, we can express (\ref{eq: L2 score equivalent renormalization physical})
in terms of the measures $\lambda$ and $\mu_{{\bf y},w}$ as follows,
\begin{align*}
{\color{teal}}S^{'}(P_{\mathbf{X}},\mathbf{y};w,h)= & \intop_{\mathbb{R}^{d}}\left(f_{\mathbf{X}}*w\right)^{2} d\lambda(\mathbf{z})  -  2\intop_{\mathbb{R}^{d}}\left(f_{\mathbf{X}}*w\right)(\mathbf{z})\,d\mu_{{\bf y},w}(\mathbf{z}).
\end{align*}
Because a nonnegative $w$ implies $f_{\mathbf{Y}}*w$ is also nonnegative,
we apply the layer cake representation to each term in the above expression
to obtain the following, 
\begin{align}
\frac{1}{2}S^{'}(P_{\mathbf{X}},\mathbf{y};w,h) & =\intop_{0}^{\infty}\Bigg(\alpha{\lambda}\Big\{ L(f_{\mathbf{X}}*w;\alpha)\Big\}-\mu_{\yy,w}\Big\{ L(f_{\mathbf{X}}*w;\alpha)\Big\}\Bigg)\,d\alpha.\label{eq:S equivalent_no last term-1}
\end{align}
In the integrand on the right-hand side of (\ref{eq:S equivalent_no last term-1}),
$L(f_{\mathbf{X}}*w;\alpha)$ is viewed as an estimate of $L(f_{\mathbf{Y}}*w;\alpha)$.
Replacing $L(f_{\mathbf{X}}*w;\alpha)$ by any Borel set $A$ leads
to (\ref{eq:scoring functions for level sets}).
}

(ii) Proof of Theorem~\ref{thm:scoring function for level set}. We need to prove that for each number $\alpha$ and each Borel set
$A\subseteq\mathbb{R}^{d}$,
\begin{equation}
\Delta':=\E_{P_{\mathbf{Y}}}\Big[{\rm (S^{'})^{\Gamma}}\big(A,\mathbf{\bullet};w,h,\alpha\big)-{\rm (S^{'})^{\Gamma}}\big(L(f_{\mathbf{Y}}*w;\alpha),\mathbf{\bullet};w,h,\alpha\big)\Big]\geq0.\label{Delta'}
\end{equation}
For notational convenience, for any $P_{\mathbf{X}}\in\mathcal{V}^{d}$
we write $\vartheta(\mathbf{y},\mathbf{z}):=w(\mathbf{z}-\mathbf{y})$,
and $\Xi[\YY](\mathbf{z}):=\left(f_{\mathbf{Y}}*w\right)(\mathbf{z})$.
We can express (\ref{Delta'}) by $\Delta'=\E_{P_{\YY}}[\Delta'']$,
where 
\begin{align*}
\Delta'' & = \underbrace{\alpha\Big(\lambda\left\{ A\right\} -\lambda\left\{ \Xi[\YY]>\alpha\right\} \Big)}_{ B_{1} } + \underbrace{\mu_{\mathbf{y},w}\big\{\Xi[\YY]>\alpha\big\}-\mu_{\mathbf{y},w}\big\{ A\big\} }_{ B_{2} }.
\end{align*}

Let us partition $\R^{d}$ as $\R^{d}=\Sigma_{++}\sqcup\Sigma_{+-}\sqcup\Sigma_{-+}\sqcup\Sigma_{--}$, where $\sqcup$ is the disjoint union, and 
\begin{align*}
 & \Sigma_{++}:=A\cap\big\{\Xi[\YY]>\alpha\big\},\qquad\Sigma_{+-}:=A\cap\big\{\Xi[\YY]\leq\alpha\big\},\\
 & \Sigma_{-+}:=A^{c}\cap\big\{\Xi[\YY]>\alpha\big\},\qquad\Sigma_{--}:=A^{c}\cap\big\{\Xi[\YY]\leq\alpha\big\},
\end{align*}
where $A^{c}$ represents the complement set of $A$. By construction,
each of these four sets is Borel measurable. Using the definitions of $d\mu_{\mathbf{y},w}$,
$d\lambda$, and Fubini's theorem,
we get 
\begin{align*}
 & B_{1}=\alpha\intop_{\mathbb{R}^{d}}\big[\I_{\Sigma_{+-}}(\zz)-\I_{\Sigma_{-+}}(\zz)\big]\,h({\bf z})\,d{\bf z},\\
 & B_{2}=\intop_{\mathbb{R}^{d}}\vartheta(\mathbf{y},\mathbf{z})\big[\1_{\Sigma_{-+}}(\mathbf{z})-\1_{\Sigma_{+-}}(\mathbf{z})\big]\,h({\bf z})\,d{\bf z},\\
 & \Delta''=\intop_{\mathbb{R}^{d}}\bigg\{\Big(\alpha-\vartheta(\mathbf{y},\mathbf{z})\Big)\big[\1_{\Sigma_{+-}}-\1_{\Sigma_{-+}}\big]\bigg\} h({\bf z})\,d{\bf z},
\end{align*}
where $\I_{B}(\mathbf{\bullet})$ denotes the indicator function on
the set $B$. Notice that $\E_{P_{\mathbf{Y}}}\left[\vartheta(\bullet,\zz)\right]=\Xi[\YY](\zz)$,
which implies that
\begin{align}
\Delta'=\E_{P_{\YY}}[\Delta'']=\intop_{\mathbb{R}^{d}}\bigg\{\Big(\alpha-\Xi[\YY](\zz)\Big)\Big[\1_{\Sigma_{+-}}(\zz)-\1_{\Sigma_{-+}}(\zz)\Big]\bigg\} h({\bf z})\,d{\bf z}.\label{final--proof 2}
\end{align}
On $\Sigma_{+-}$ one has $\1_{\Sigma_{+-}}-\1_{\Sigma_{-+}}=1$ and
$\alpha-\Xi[\YY]\geq0$, and on $\Sigma_{-+}$, $\1_{\Sigma_{+-}}-\1_{\Sigma_{-+}}=-1$
and $\alpha-\Xi[\YY]\leq0$. Therefore, the integrand of $\Delta'$
is pointwise non-negative, so $\Delta'\geq0$.

Regarding the strict consistency, we adopt the approach considered
by \citet{fissler2019forecast}. If $A$ and $L(f_{\mathbf{Y}}*w;\alpha)$
differ by a set $B$ that is not $\mathcal{L}^{d}$-null, then
\begin{equation}
\Delta'=\intop_{B}\bigg\{\Big(\alpha-\Xi[\YY](\zz)\Big)\Big[\1_{\Sigma_{+-}}(\zz)-\1_{\Sigma_{-+}}(\zz)\Big]\bigg\} h({\bf z})\,d{\bf z}.\label{eq: Final-proof 3}
\end{equation}
However, if $L(f_{\mathbf{Y}}*w;\alpha)$ coincides with the closure
of $\big\{ f_{{\bf Y}}\ast w>\alpha\big\}$, the set $\{f_{{\bf Y}}\ast w=\alpha\}$
is $\mathcal{L}^{d}$-null. It follows that the integrand on the right-hand
side of (\ref{eq: Final-proof 3}) is strictly positive on $B$ modulo
an $\mathcal{L}^{d}$-null set. But by Assumption~\ref{assumption: one-1}(c), one has $h>0$ $\mathcal{L}^{d}$-a.e.; so $\Delta'$ is strictly
positive. This proves the strict consistency of $\left(S'\right)^{\Gamma}$. \Halmos

\endproof

\end{APPENDIX}

\begin{APPENDIX}{Discussion of LPMs in Section \ref{subsec: Scoring Function for CDF level sets}}

\textcolor{\mycolor}{We first show that LPMs uniquely characterize distributions. To see this, we can take the $k^{th}$ partial derivatives
with respect to all $z_{j}$,
\begin{align*}
\left[\frac{\partial^{k}}{\partial z_{1}^{k}}\cdots\frac{\partial^{k}}{\partial z_{d}^{k}}(f_{\mathbf{X}}\ast u^{{*k}})\right](\zz) & =\int_{-\infty}^{z_{d}}\cdots\int_{-\infty}^{z_{1}}\frac{\partial^{k}}{\partial s_{1}^{k}}\cdots\frac{\partial^{k}}{\partial s_{d}^{k}}u^{{*k}}(\mathbf{s})f_{\mathbf{X}}(s_{1},\ldots,s_{d})\,ds_{1}\cdots ds_{d}\\
 & =\int_{-\infty}^{z_{d}}\cdots\int_{-\infty}^{z_{1}}f_{\mathbf{X}}(s_{1},\ldots,s_{d})\,ds_{1}\cdots ds_{d}=f_{\mathbf{X}}\ast u(\zz)=F_{\mathbf{\mathbf{X}}}(\zz).
\end{align*}
Therefore, if $\text{LPMS}^{'}(P_{\mathbf{X}},\mathbf{y};k,h)=\text{LPMS}^{'}(P_{\mathbf{Y}},\mathbf{y};k,h)$,
then their CDFs must coincide, which means $P_{\mathbf{X}}=P_{\mathbf{Y}}$.
We also note that when $k=0$, ${\rm LPM}_{{\bf X},0}=F_{\mathbf{X}}$
is simply the CDF of $P_{\mathbf{X}}$, and in this case $\text{LPMS}$
and $\text{LPMS}^{'}$ are simply the MCRPS and $\text{MCRPS}^{'}$
discussed in Section~\ref{subsec: MCRPS}.}

\textcolor{\mycolor}{In addition to the three scores in Section \ref{subsec: Examples of Scoring Functions},
the following two conditions can also guarantee that $f_{\mathbf{Y}}*w$
uniquely characterizes $P_{{\bf Y}}$: }

\textcolor{\mycolor}{(i) if $w$ and every distribution of $\mathcal{V}^{d}$
are compactly supported: it follows from J.-L. Lions' generalisation
of the Titchmarsh convolution theorem (\citealt{lions1951supports})
that $f_{{\bf X}}=f_{{\bf Y}}$ $\mathcal{L}^{d}$-a.e..}

\textcolor{\mycolor}{(ii) Let $w$ and/or its certain weak derivatives
have well-defined Fourier transforms (denoted with a hat) that are nonzero $\mathcal{L}^{d}$-a.e.. First, if $\hat{w}$ exists and is nonzero $\mathcal{L}^{d}$-a.e.,
then $f_{\mathbf{X}}*w=f_{\mathbf{Y}}*w$ $\mathcal{L}^{d}$-a.e.
implies $f_{\mathbf{X}}=f_{\mathbf{Y}}$ $\mathcal{L}^{d}$-a.e..
Next, let $v$ be any given weak derivative of $w$ and $\hat{v}\neq0$
$\mathcal{L}^{d}$-a.e.. Then, the corresponding derivatives of $f_{X}*w$
and $f_{Y}*w$ are precisely $f_{X}*v$ and $f_{Y}*v$, respectively,
thanks to the fundamental theorem of calculus. Finally, for $f_{\mathbf{X}}*w=f_{\mathbf{Y}}*w$
(or, using the fact that convolution intertwines with multiplication,
we obtain $\hat{f}_{\mathbf{X}}\hat{w}=\hat{f}_{\mathbf{Y}}\hat{w}$.
Because $\hat{w}\neq0$ $\mathcal{L}^{d}$-a.e., we conclude $\hat{f}_{\mathbf{X}}=\hat{f}_{\mathbf{Y}}$
$\mathcal{L}^{d}$-a.e., which implies $P_{\mathbf{X}}=P_{\mathbf{Y}}$.
The same argument holds for $v$.}

\end{APPENDIX}

\begin{APPENDIX}{Discussion on Remark \ref{Remark: Fourier side L2 score}}

For a function $g:\R^{d}\map\CC$, its Fourier transform $\hat{g}\equiv\mathcal{F}(g):\R^{d}\map\CC$
is given by $\hat{g}(\mathbf{t}):=\int_{\R^{d}}e^{2\pi i\mathbf{z}^{T}\mathbf{t}}g(\mathbf{z})\,d\mathbf{z},$
where $\mathbf{t}\in\R^{d}$ and $i$ is the imaginary unit. When
$g$ is a PDF of a distribution, $\hat{g}$ is called the characteristic
function of the distribution. The definition of the Fourier transform
can be extended to certain generalized functions including Dirac delta
masses.

Assume that the Plancherel identity holds. Then, using the fact that
convolution intertwines with multiplication, we can obtain a Fourier
representation for the $L^{2}$ score in (\ref{eq: L2 Score physical}),
\begin{align}
S(P_{\mathbf{X}},\mathbf{y};\hat{w},\hat{h}) & =\intop_{\mathbb{R}^{d}}\Big|\big(\hat{f}_{\mathbf{X}}\hat{w}-\hat{\delta}_{\mathbf{y}}\hat{w}\big)*\hat{h}\Big|^{2}(\mathbf{t})\,d\mathbf{t},\label{eq: L^2 score FOurier}
\end{align}
where $\hat{w},\hat{h},\hat{f}_{\mathbf{X}},\hat{\delta}_{\mathbf{y}}$
denote the Fourier transforms of $w,h,f_{\mathbf{X}},\delta_{\mathbf{y}}$,
respectively.

Note that (\ref{eq: L^2 score FOurier}) can remain valid even if
the Plancherel identity does not hold, provided that $S(P_{\mathbf{X}},\mathbf{y};\hat{w},\hat{h})$
is finite for all $\mathbf{y}$ and that $S(P_{\mathbf{X}},\mathbf{\bullet};\hat{w},\hat{h})f_{\mathbf{Y}}(\mathbf{\bullet})$
is integrable for any $P_{\mathbf{Y}}\in\mathcal{V}^{d}$. In general,
there is no symmetry between (\ref{eq: L2 Score physical}) and (\ref{eq: L^2 score FOurier}).
In Examples~\ref{subsec: MCRPS} and \ref{subsec:Lower-Partial-Moments scoring function},
$|w|$ is not integrable hence does not have Fourier transform, hence
(\ref{eq: L^2 score FOurier}) cannot be defined; conversely, if we
consider $\hat{w}$ with no well-defined inverse Fourier transform,
then one does not have (\ref{eq: L2 Score physical}).

If we consider $\hat{w}_{\text{ES}}(\mathbf{t})=\big(\lVert\mathbf{t}\lVert^{\frac{d+1}{2}}\big)^{-1}/2$
and $\hat{h}(\mathbf{t})=\delta_{\mathbf{0}}(\mathbf{t})$, (\ref{eq: L^2 score FOurier})
leads to the well-known energy score (see, for example, \citealt{baringhaus2004new,gneiting2007strictly}).
For the level set score, as the inverse Fourier transform $w=\mathcal{F}^{-1}(\hat{w})$
does not exists, we cannot directly apply Theorem~\ref{thm:scoring function for level set}.
Observing that $\hat{w}_{\text{ES}}(\mathbf{t})$
is a radial function, we can consider an alternative ``projective''
approach, where we take the Fourier transform of $\hat{w}$ and apply
Theorem~\ref{thm:scoring function for level set} ``projectively''
along different directions over the unit sphere (see, for example,
\citealt{baringhaus2004new}). This leads to the scoring function
for the so-called projective quantile (\citealt{kong2012quantile}).
However, this approach is essentially straightforward, as it only
relies on the decomposition of the multivariate distributions into
univariate projections. This is not the focus of the paper, and so
we do not discuss this further.

\end{APPENDIX}

\end{document}